\documentclass[11pt, a4paper]{amsart}
\usepackage{amsmath, amsthm, amssymb}
\usepackage{txfonts, mathrsfs}
\usepackage[usenames]{color}
\usepackage{psfrag}
\usepackage{graphicx}
\usepackage{todonotes}

\usepackage{comment}

\numberwithin{equation}{section}

\newtheorem{thm}{Theorem}[section]
\newcommand{\bt}{\begin{thm}}
\newcommand{\et}{\end{thm}}

\newtheorem{cor}[thm]{Corollary}
\newcommand{\bc}{\begin{cor}}
\newcommand{\ec}{\end{cor}}

\newtheorem{lem}[thm]{Lemma}
\newcommand{\bl}{\begin{lem}}
\newcommand{\el}{\end{lem}}

\newtheorem{prop}[thm]{Proposition}
\newcommand{\bp}{\begin{prop}}
\newcommand{\ep}{\end{prop}}

\newtheorem{defn}[thm]{Definition}
\newcommand{\bd}{\begin{defn}}
\newcommand{\ed}{\end{defn}}

\newtheorem{rmrk}[thm]{Remark}
\newcommand{\br}{\begin{rmrk}}
\newcommand{\er}{\end{rmrk}}

\newtheorem{quest}[thm]{Question}
\newcommand{\bq}{\begin{quest}}
\newcommand{\eq}{\end{quest}}

\newtheorem{example}[thm]{Example}

\newcommand{\N}{\mathbb{N}}

\newdimen\vintkern\vintkern12pt
\def\vint{-\kern-\vintkern\int}

\newcommand{\hm}{{\mathcal H}}

\newcommand{\diam}{\operatorname{diam}}
\newcommand{\trace}{\operatorname{tr}}
\newcommand{\length}{\ell}
\newcommand{\Area}{\operatorname{Area}}

\newcommand{\md}{\operatorname{md}}

\newcommand{\sys}{\operatorname{sys}}

\newcommand{\jac}{{\mathbf J}}
\newcommand{\ap}{\operatorname{ap}}
\newcommand{\apmd}{\ap\md}
\newcommand{\inter}{\operatorname{int}}

\newcommand{\R}{\mathbb{R}}

\begin{document}
\pagebreak
\bibliographystyle{plain}


\title{Area minimizing surfaces in homotopy classes in metric spaces}


\subjclass[2010]{49Q05 (53A10, 53C23)}

\author{Elefterios Soultanis}

\address
  {Department of Mathematics\\ University of Fribourg\\ Chemin du Mus\'ee 23\\ 1700 Fribourg, Switzerland}
\email{elefterios.soultanis@gmail.com}

\author{Stefan Wenger}

\address
  {Department of Mathematics\\ University of Fribourg\\ Chemin du Mus\'ee 23\\ 1700 Fribourg, Switzerland}
\email{stefan.wenger@unifr.ch}

\date{\today}

\thanks{Research supported by Swiss National Science Foundation Grants 165848 and 182423}

\begin{abstract}
 We introduce and study a notion of relative $1$--homotopy type for Sobolev maps from a surface to a metric space spanning a given collection of Jordan curves. We use this to establish the existence and local H\"older regularity of area minimizing surfaces in a given relative $1$--homotopy class in proper geodesic metric spaces admitting a local quadratic isoperimetric inequality. If the underlying space has trivial second homotopy group then relatively $1$--homotopic maps are relatively homotopic. We also obtain an analog for closed surfaces in a given $1$--homotopy class. Our theorems generalize and strengthen results of Lemaire, Jost, Schoen-Yau, and Sacks-Uhlenbeck. 
\end{abstract}

\maketitle

\renewcommand{\theequation}{\arabic{section}.\arabic{equation}}
\pagenumbering{arabic}

\section{Introduction}\label{sec:Intro}

\subsection{Background}
Let $M$ be a $2$--dimensional surface with boundary. A map from $M$ to a Riemannian manifold $N$ is said to span a given collection $\Gamma\subset N$ of Jordan curves if its restriction to $\partial M$ is a weakly monotone parametrization of $\Gamma$. Consider the problem of finding a weakly conformal map of minimal area among maps spanning $\Gamma$. When $M$ is a disc this amounts to the classical Problem of Plateau, with first general solutions going back to \cite{Dou31, Rad30, Cou37} for $N=\R^n$ and to \cite{Mor48} for homogeneously regular Riemannian manifolds $N$. When $M$ is a surface of higher topological type, possibly with several boundary components, the problem is known as the Plateau-Douglas problem. It was first considered in \cite{Dou39, Shi39, Cou40} with different non-degeneracy conditions; complete modern solutions appeared in \cite{Jos85,TT88}. 

One may further ask whether it is possible to find a weakly conformal map of minimal area spanning $\Gamma$ in a fixed relative homotopy class. In general, such maps need not exist, see \cite{Jost84, Lem78}. However, Lemaire \cite{Lem82} showed the existence of an area minimizer in a fixed relative homotopy class under the assumption that $N$ has trivial second homotopy group, while Jost \cite{Jos85} proved the existence of an area minimizer inducing the same action on fundamental groups as a given map. Schoen-Yau \cite{SY79} and Sacks-Uhlenbeck \cite{SU82} considered the related problem of finding a mapping of minimal area from a closed (i.e. compact and without boundary) surface $M$  to $N$ inducing the same action on fundamental groups as a given map. Finally, White \cite{White88} introduced the notion of $d$--homotopy type for Sobolev maps from a closed manifold of any dimension to a Riemannian manifold and proved the existence of mappings of minimal energy in a given $d$--homotopy class for suitable integers $d$.

Recently, the classical Plateau and the Plateau-Douglas problems have been solved in metric spaces of various generality in \cite{Nik79, Jost94, MZ10, OvdM14, LW15-Plateau, Cre-Plateau-sing} and \cite{FW-Plateau-Douglas,  Creutz-Fitzi}, respectively. 
In the present article we strengthen the results of Lemaire \cite{Lem82} and Jost \cite{Jos85} mentioned above and generalize them to the setting of proper geodesic metric spaces admitting a local quadratic isoperimetric inequality. For this purpose, we introduce and study a notion of $1$--homotopy classes of Sobolev maps relative to a given collection of Jordan curves. Our notion is akin to $d$--homotopy of Sobolev maps defined on a closed manifold introduced by White \cite{White88} and studied in \cite{White88, HL03}. It provides better control than the induced action on the fundamental group. We then solve the Plateau-Douglas problem in relative $1$--homotopy classes and show that solutions are locally H\"older continuous and conformal in a weak metric sense.  If the underlying space has trivial second homotopy group then relatively $1$--homotopic maps are relatively homotopic. To our knowledge, our results are already partially new for Riemannian manifolds. We further obtain an analog for closed surfaces, generalizing the results in \cite{SY79, SU82} mentioned above.

\bd\label{def:qii}
A complete metric space $X$ is said to admit a local quadratic isoperimetric inequality if there exist $C, l_0>0$ such that every Lipschitz curve $c\colon S^1\to X$ of length $\length(c)\leq l_0$ is the trace of a Sobolev map $u\in W^{1,2}(\mathbb D, X)$ with $$\Area(u)\leq C\cdot \length(c)^2.$$
\ed

For the notions related to Sobolev maps we refer to Section~\ref{sec:prelims}.
The class of spaces admitting a local quadratic isoperimetric inequality contains all homogeneously regular Riemannian manifolds \cite{Mor48}, compact Lipschitz manifolds, complete $\rm{CAT}(\kappa)$--spaces, compact Alexandrov spaces, some sub-Riemannian manifolds, and many more spaces, cf. \cite[Section 8]{LW15-Plateau}.

\subsection{Relative $1$--Homotopy classes of Sobolev maps}

Let $\Gamma\subset X$ be the disjoint union of $k\ge 1$ rectifiable Jordan curves in a proper geodesic metric space $X$ admitting a local quadratic isoperimetric inequality. Let $M$ be a smooth compact oriented surface with $k$ boundary components, and let $g$ be an auxiliary Riemannian metric on $M$. We denote by $[\Gamma]$ the family of weakly monotone  parametrizations of $\Gamma$, i.e. uniform limits of homeomorphisms $\partial M\to \Gamma$, and by $\Lambda(M,\Gamma,X)$ the family of Sobolev maps $u\in W^{1,2}(M, X)$ such that the trace $\trace(u)$ has a continuous representative in $[\Gamma]$. Let $h\colon K\to M$ be a $C^1$--smooth triangulation of $M$, and $\varrho\colon K^1\to X$ a continuous map such that $\varrho|_{\partial K}\in [\Gamma]$, where $K^1$ denotes the $1$--skeleton of $K$ and $\partial K\subset K^1$ is the subset of $K$ homeomorphic to $\partial M$. The  homotopy class of $\varrho$ relative to $\Gamma$ is the family $$[\varrho]_\Gamma:= \{\varrho'\colon K^1\to X \mid\  \varrho'\textrm{ continuous },\ \varrho'|_{\partial K}\in[\Gamma],\ \varrho\sim\varrho'\textrm{ rel }\Gamma\},$$ where $\varrho$ and $\varrho'$ are said to be homotopic relative to $\Gamma$, denoted $\varrho\sim \varrho'$ rel $\Gamma$, if there exists a homotopy $F\colon K^1\times[0,1]\to X$ from $\varrho$ to $\varrho'$ with $F(\cdot, t)|_{\partial K}\in [\Gamma]$ for every $t$.

The $1$--homotopy class $u_{\#, 1}[h]$ relative to $\Gamma$ of an element $u\in \Lambda(M,\Gamma, X)$ will be defined in Section~\ref{sec:1-homot}. In the following theorem we summarize its most important properties. These could in fact be used to give an equivalent definition of $u_{\#, 1}[h]$, see the remark after the theorem. 

\bt\label{thm:intro-properties-1-hom-class-Sobolev}
 Every $u\in \Lambda(M, \Gamma, X)$ has a well-defined relative homotopy class $u_{\#, 1}[h]$ of continuous maps from $K^1$ to $X$ whose restriction to $\partial K$ is in $[\Gamma]$. It satisfies:
 \begin{enumerate}
  \item If $u$ has a representative $\bar{u}$ which is continuous on the whole of $M$ then $$u_{\#,1}[h] = [\bar{u}\circ h|_{K^1}]_\Gamma.$$
    \item If $u,v\in \Lambda(M,\Gamma, X)$ satisfy $u_{\#,1}[h] = v_{\#, 1}[h]$ then, for every triangulation $\tilde{h}\colon \tilde{K}\to M$ of $M$, we have $$u_{\#,1}[\tilde{h}] = v_{\#,1}[\tilde{h}].$$
  \item For every $L>0$ there exists $\varepsilon>0$ such that if $u,v\in\Lambda(M, \Gamma, X)$ induce the same orientation on $\Gamma$,  and $$d_{L^2}(u,v)\leq \varepsilon,\quad \max\left\{E_+^2(u,g), E_+^2(v,g)\right\}\leq L,$$  then $u_{\#,1}[h] = v_{\#,1}[h].$
 \end{enumerate}
\et

Here, $E_+^2(u,g)$ denotes the Reshetnyak energy of $u$ with respect to $g$, see Section~\ref{sec:prelims}. Maps in $\Lambda(M,\Gamma, X)$ can be approximated in the $L^2$--distance by \emph{continuous} maps in $\Lambda(M,\Gamma, X)$ with the same trace and control on the energy, see Lemma~\ref{lem:approx-cont-area-energy-bounded}. Thus properties (i) and (iii) in Theorem~\ref{thm:intro-properties-1-hom-class-Sobolev} imply that the $1$--homotopy class $u_{\#,1}[h]$ is well-defined. The argument used to prove (ii) also shows that, if $u\in \Lambda(M,\Gamma, X)$ and $\varphi\colon M\to X$ is continuous with $\varphi|_{\partial M}\in[\Gamma]$, then $u_{\#,1}[h] = [\varphi\circ h|_{K^1}]_\Gamma$ holds for one triangulation $h$ if and only if it holds for every triangulation. In this case we say that $u\in\Lambda(M,\Gamma, X)$ is $1$--homotopic to $\varphi$ relative to $\Gamma$, denoted by $u\sim_1 \varphi$ rel $\Gamma$.

\subsection{Homotopic Plateau-Douglas problem}
Let $\Gamma$, $X$ be as above, and let $M$ be a smooth compact oriented and \emph{connected} surface with $k\ge 1$ boundary components. Given a continuous map $\varphi\colon M\to X$ with $\varphi|_{\partial M}\in[\Gamma]$, set 
$$a(M, \varphi, X):= \inf\{\Area(u): \text{$u\in\Lambda(M, \Gamma, X)$, $u\sim_1\varphi$ rel $\Gamma$}\},$$
where $\inf\varnothing=\infty$ by convention. Moreover, set $a^*(M, \varphi, X):= \inf a(M^*, \varphi^*, X)$,
where the infimum is taken over all \emph{primary reductions} of $(M,\varphi)$, that is, pairs $(M^*, \varphi^*)$ consisting of
\begin{enumerate}
\item[(i)] a smooth surface $M^*$ obtained from $M$ by cutting $M$ along a smooth closed simple non-contractible curve $\alpha$ in the interior of $M$ and gluing smooth discs to the two new boundary components;
\item[(ii)] a continuous map $\varphi^*\colon M^*\to X$ which agrees with $\varphi$ on $M\setminus \alpha$. 
\end{enumerate}
We say that $\varphi$ satisfies the \emph{homotopic Douglas condition} if 
\begin{equation}\label{eq:Douglas-condition}
 a(M,\varphi,X) < a^*(M,\varphi,X).
\end{equation}
As an illustration, if the induced homomorphism $\varphi_*\colon \pi_1(M) \to \pi_1(X)$ of fundamental groups is injective then $\varphi$ satisfies the homotopic Douglas condition \eqref{eq:Douglas-condition} and, in par\-ti\-cu\-lar, $a(M,\varphi,X)<\infty$, see Proposition~\ref{prop:induced-homom-Douglas}. In the statement below, we fix $\Gamma,\ X$, and $M$ as above, and let $\varphi\colon M\to X$ be a continuous map with $\varphi|_{\partial M}\in [\Gamma]$.

\bt\label{thm:Plateau-Douglas-homot-intro}
 If $\varphi$ satisfies the homotopic Douglas condition \eqref{eq:Douglas-condition} then:
 \begin{enumerate}
  \item There exist $u\in\Lambda(M, \Gamma, X)$ and a Riemannian metric $g$ on $M$ such that $u$ is $1$--homotopic to $\varphi$ relative to $\Gamma$, $u$ is infinitesimally isotropic with respect to $g$, and $\Area(u) = a(M, \varphi, X)$. 
  \item Any such $u$ has a representative $\bar{u}$ which is locally H\"older continuous in the interior of $M$ and extends continuously to the boundary $\partial M$.
  \item If $X$ has trivial second homotopy group then $\bar{u}$ is homotopic to $\varphi$ relative to $\Gamma$.
 \end{enumerate}
\et

Moreover, the metric $g$ can be chosen such that it has constant curvature $-1$, $0$, or $1$ and $\partial M$ is geodesic. See Section~\ref{sec:prelims} for the definition of infinitesimal isotropy, which is a metric variant of weak conformality. Here, $\bar{u}$ and $\varphi$ are called homotopic relative to $\Gamma$ if they are homotopic through a family of maps whose restriction to $\partial M$ is in $[\Gamma]$. We remark that homotopy classes (relative to $\Gamma$) need not contain continuous, infinitesimally isotropic area minimizers if $\pi_2(X)\ne\varnothing$, compare \cite[Chapter 5]{Jost84}. 

Theorem~\ref{thm:Plateau-Douglas-homot-intro} generalizes and strengthens \cite[Theorem 2.2]{Jos85} and \cite[Theorem 1.7]{Lem82}, see also \cite[Theorem 5.1]{Jost94} for a homotopic variant of the Dirichlet problem in metric spaces. We remark that control on the relative $1$--homotopy class is, in general, strictly stronger than the control on the action on fundamental groups in \cite{Jos85}, see Example~\ref{ex:1-hom-stronger-action-fundgrp}. An analog of Theorem~\ref{thm:Plateau-Douglas-homot-intro} for closed surfaces, generalizing results in \cite{SY79,SU82}, will be discussed in Section~\ref{sec:sol}.

We remark that the local quadratic isoperimetric inequality is crucial to the stability statement (iii) in Theorem~\ref{thm:intro-properties-1-hom-class-Sobolev}. Example~\ref{ex:stabilityfail} exhibits a space where the stability of $1$--homotopy classes from closed surfaces fails. Compare with \cite{Creutz-Fitzi}, where the Plateau-Douglas problem was recently solved in spaces without a local quadratic isoperimetric inequality.

\subsection{Outline} 

The idea for defining the relative $1$--homotopy type of a map $u\in\Lambda(M,\Gamma,X)$ is, like in \cite{White88}, to consider small perturbations of $C^1$--smooth triangulations of $M$ in such a way that the restriction of $u$ to the $1$--skeleton of a ``generic'' perturbed triangulation is essentially continuous. In Section~\ref{sec:admissible-deformations}, we introduce \emph{admissible deformations} on $M$ which accomplish this and prove that the relative homotopy class of such restrictions is essentially independent of the perturbation, see Theorem~\ref{thm:homotopic-1-skeleton}. This crucially uses the local quadratic isoperimetric inequality. 

In Section~\ref{sec:1-homot} we show that the way we perturb a given triangulation does not affect the relative homotopy type of the restrictions to generic $1$--skeleta. Together with a continuous approximation of Sobolev maps (see Lemma~\ref{lem:approx-cont-area-energy-bounded}) and the results of Section~\ref{sec:admissible-deformations}, this leads to a well-defined notion of relative $1$--homotopy class for Sobolev maps, which is moreover independent of the chosen triangulation. The main results in Section~\ref{sec:1-homot} are Theorems~\ref{thm:rel-hom-indep-wiggling} and \ref{thm:stability-1-homotopic} from which Theorem~\ref{thm:intro-properties-1-hom-class-Sobolev} will follow. As already mentioned, our notion of relative $1$--homotopy class is related to the $d$--homotopy type, studied primarily for Sobolev maps defined on closed manifolds in \cite{White88, HL03, HL05}. While these articles also discuss the case of manifolds with boundary, Sobolev maps in their setting are required to have a fixed Lipschitz trace. This is suitable for solving the Dirichlet problem in $d$--homotopy classes but cannot be applied to the Plateau--Douglas problem since it is not possible to control the boundary behaviour of elements of $\Lambda(M,\Gamma, X)$.

In Sections~\ref{sec:homot-Douglas-cond} and \ref{sec:sol} we use an approach analogous to that in \cite{FW-Plateau-Douglas} in order to solve the homotopic Plateau-Douglas problem. Unlike in \cite{FW-Plateau-Douglas}, we need to control the relative $1$--homotopy type of the primary reductions appearing in the proofs of Propositions~\ref{prop:equi-cont} and \ref{prop:lower-bound-rel-systole}. Lemma~\ref{lem:1-hom-reduction} provides the necessary technical tool for this. We furthermore provide a simple sufficient condition (see Proposition~\ref{prop:induced-homom-Douglas}) that ensures the homotopic Douglas condition \eqref{eq:Douglas-condition} is satisfied. Section~\ref{sec:sol} is devoted to the proof of Theorem~\ref{thm:Plateau-Douglas-homot-intro}. We present and prove Theorem~\ref{thm:area-min-hom-class-without-bdry}, which is an analog of Theorem~\ref{thm:Plateau-Douglas-homot-intro} for closed surfaces.

\section{Preliminaries}\label{sec:prelims}

\subsection{Terminology}

A \emph{surface}, in this work, refers to a smooth compact oriented surface with (possibly) non-empty boundary, and a \emph{closed surface} is a surface with empty boundary. We denote by $\partial M$ and  ${\rm int}(M)=M\setminus \partial M$ the boundary and interior of a surface $M$, respectively. The Euler characteristic of a connected surface $M$ satisfies $\chi(M)=2-2p-k$, where $k\ge 0$ is the number of components of $\partial M$, and $p$ is the genus of the closed surface obtained by gluing a disc along every boundary component of $M$.

For a metric space $X$ and $m\ge 0$, we denote by $\hm_X^m$ the Hausdorff $m$--measure on $X$. If $X$ is a manifold equipped with a Riemannian metric $g$, we denote $\hm_g^m=\hm_X^m$. The Lebesgue measure of a subset $A\subset \R^m$ is denoted by $|A|$.

\subsection{Triangulations}

A triangulation of a surface $M$ is a homeomorphism $h\colon K\to M$ from a cell-complex $K$, equipped with the length metric which restricts to the Euclidean metric on every cell $\Delta$ of $K$. We additionally assume throughout the paper that triangulations are $C^1$--diffeomorphisms, i.e. $h|_\Delta$ is a $C^1$--diffeomorphism onto its image for any cell $\Delta$ of $K$ (cells are closed by definition). The $j$--skeleton $K^j$ of $K$ is the union of the cells of $K$ with dimension $\le j$, and $\partial K\subset K^1$ is the subset of $K$ homeomorphic to $\partial M$.

\subsection{Semi-norms}

The energy of a semi-norm $s$ on (Euclidean) $\R^2$ is defined by  $$\mathbf{I}_+^2(s):= \max\{s(v)^2: v\in\R^2, |v|=1\}.$$
The jacobian of a norm $s$ on $\R^2$ is the unique number $\jac(s)$ such that $$\hm^2_{(\R^2, s)}(A) = \jac(s) \cdot |A|$$ for some and thus every subset $A\subset\R^2$ with $|A|>0$. For a degenerate semi-norm $s$ we set $\jac(s):= 0$. Notice that we always have $\jac(s)\leq \mathbf{I}_+^2(s)$. A semi-norm $s$ on $\R^2$ is called isotropic if $s=0$ or if $s$ is a norm and the ellipse of maximal area contained in $\{v\in\R^2: s(v)\leq 1\}$ is a round Euclidean ball.

\subsection{Sobolev maps with metric targets}

Let $(X, d)$ be a complete metric space and let $M$ be a smooth compact $m$--dimen\-sional manifold, possibly with non-empty boundary. 
Fix a Riemannian metric $g$ on $M$ and let $\Omega\subset M$ be open and bounded. 

Denote by $L^2(\Omega, X)$ the collection of measurable and essentially separably valued maps $u\colon \Omega\to X$ such that for some and thus every $x\in X$ the function $u_x(z):= d(x, u(z))$ belongs to the classical space $L^2(\Omega)$. For $u,v\in L^2(\Omega, X)$ we define $$d_{L^2}(u,v):= \left(\int_\Omega d^2(u(z), v(z))\,d\hm_g^m(z)\right)^{\frac{1}{2}},$$
and we say that a sequence $(u_n)\subset L^2(\Omega, X)$ converges in $L^2(\Omega, X)$ to $u\in L^2(\Omega, X)$ if $d_{L^2}(u_n,u)\to 0$ as $n\to\infty$. The following definition is due to Reshetnyak \cite{Res97,Res06}.

\bd
 A map $u\in L^2(\Omega, X)$ belongs to the Sobolev space $W^{1,2}(\Omega, X)$ if there exists $h\in L^2(\Omega)$ such that $u_x\in W^{1,2}({\rm int}(M))$ and $|\nabla u_x|_g\leq h$ almost everywhere on $\Omega$, for every $x\in X$.
\ed

Several other notions of Sobolev spaces exist in the literature and we refer the reader to \cite[Chapter 10]{HKST15} for an overview of some of them. We will use in particular \emph{Newton-Sobolev} spaces which are equivalent to $W^{1,2}(\Omega,X)$ if $\Omega$ is a bounded Lipschitz domain, see Proposition~\ref{prop:Newton-Sobolev-rep} for a precise statement.

If $u\in W^{1,2}(\Omega, X)$ then for almost every $z\in \Omega$ there exists a unique semi-norm $\apmd u_z$ on $T_zM$ such that $$\ap \lim_{v\to 0} \frac{d(u(\exp_z(v), u(z)) - \apmd u_z(v)}{|v|_g}=0,$$
where $\ap\lim$ is the approximate limit, see e.g. \cite{Kar07}. Next, we specialize to the case that $M$ has dimension $m=2$. We define the notions of energy, jacobian and isotropy of a semi-norm on $(T_zM, g(z))$ by identifying it with $(\R^2,|\cdot|)$ via a linear isometry.

\bd
 Let $u\in W^{1,2}(\Omega, X)$. The Reshetnyak energy of $u$ with respect to $g$ and the parametrized (Hausdorff) area of $u$  are given, respectively, by $$E_+^2(u, g):= \int_{\Omega} \mathbf{I}_+^2(\apmd u_z)\,d\hm^2_{g}(z),\quad  \Area(u):= \int_{\Omega} \jac(\apmd u_z)\,d\hm^2_g(z).$$
\ed

We have that the parametrized area of a Sobolev map is invariant under precompositions with biLipschitz homeo\-morphisms, and thus independent of the Riemannian metric $g$. The energy $E_+^2$ is invariant only under precompositions with conformal diffeomorphisms, and thus depends on $g$. Our notation reflects these facts. 
Finally, if $u$ satisfies Lusin's property (N) then the area formula \cite{Kir94}, \cite{Kar07} for metric space valued Sobolev maps yields $$\Area(u) = \int_X\#u^{-1}(x)\,d\hm^2_X(x).$$

\bd
 A map $u\in W^{1,2}(\Omega, X)$ is called infinitesimally isotropic with respect to the Riemannian metric $g$ if for almost every $z\in \Omega$ the semi-norm $\apmd u_z$ on $(T_zM, g(z))$ is isotropic.
\ed

If $X$ is a Riemannian manifold, or more generally a space with property (ET) (cf. \cite[Definition 11.1]{LW15-Plateau}), then infinitesimal isotropy is equivalent to weak conformality, see \cite[Theorem 11.3]{LW15-Plateau}.

Next, we recall the definition of the trace of a Sobolev map. Let $\Omega \subset {\rm int}(M)$ be a Lipschitz domain. Then for every $z\in \partial \Omega$ there exist an open neighborhood $U\subset M$ and a biLip\-schitz map $\psi\colon (0,1)\times [0,1)\to M$ such that $\psi((0,1)\times (0,1)) = U\cap \Omega$ and $\psi((0,1)\times\{0\}) = U\cap \partial\Omega$. Let $u\in W^{1,2}(\Omega, X)$. For almost every $s\in (0,1)$ the map $t\mapsto u\circ\psi(s,t)$ has an absolutely continuous representative which we denote by the same expression. The trace of $u$ is defined by $$\trace(u)(\psi(s,0)):= \lim_{t\searrow 0} (u\circ\psi)(s,t)$$ for almost every $s\in(0,1)$. It can be shown (see \cite{KS93}) that the trace is independent of the choice of the map $\psi$ and defines an element of $L^2(\partial \Omega, X)$.

\bp\label{prop:good-cont-filling}
 Let $X$ be a proper metric space admitting a local quadratic isoperimetric inequality. Let $\Omega$ be a Lipschitz Jordan domain in the interior of $M$ and let $u\in W^{1,2}(\Omega, X)$ have a continuous trace. Then for every $\varepsilon>0$ there exists a continuous map $v\colon\overline{\Omega}\to X$ with $v|_{\partial \Omega} = \trace(u)$,  $v\in W^{1,2}(\Omega, X)$, and 
 \begin{align*}
 \Area(v)\leq \Area(u) + \varepsilon\cdot E_+^2(u, g),\quad
 E_+^2(v, g) \leq \left(1+\varepsilon^{-1}\right)\cdot E_+^2(u, g).
 \end{align*}
\ep

It follows, in particular, that if a closed curve $\gamma$ in $X$ is the trace of a Sobolev disc then $\gamma$ is contractible.

\begin{proof}
  By possibly doubling $M$ we may assume that $M$ has no boundary. Now, there exists a conformal diffeomorphism from a bounded open subset of $\R^2$ onto an open subset of $M$ which contains $\overline{\Omega}$. 
 Since area and energy are invariant under conformal diffeomorphisms we may assume that $\Omega$ is a bounded Lipschitz Jordan domain in $\R^2$. We write $E_+^2(u)$ for the energy of $u$.
 
Fix $\varepsilon>0$. We first show the existence of a minimizer $v\in W^{1,2}(\Omega, X)$ of $$A_\varepsilon(v):= \Area(v) + \varepsilon\cdot E_+^2(v),$$ subject to the condition $\trace(v) = \trace(u)$. For this let $(v_n)\subset W^{1,2}(\Omega, X)$ be a minimizing sequence for $A_\varepsilon$ with $\trace(v_n) = \trace(u)$ for all $n$. Then $(v_n)$ has bounded energy and thus, by \cite[Lemma 4.11]{LW15-Plateau} and \cite[Theorems 1.13 and 1.12.2]{KS93}, a subsequence converges in $L^2(\Omega, X)$ to a map $v\in W^{1,2}(\Omega, X)$ with $\trace(v) = \trace(u)$. By the lower semi-continuity of area and energy it follows that $v$ is a minimizer of $A_\varepsilon$.

Next, we claim that for every Lipschitz domain $\Omega'\subset \Omega$ and every $w\in W^{1,2}(\Omega', X)$ with $\trace(w) = \trace(v|_{\Omega'})$ we have $E_+^2(v|_{\Omega'}) \leq \left(1+\varepsilon^{-1}\right)\cdot E_+^2(w)$ and thus $v$ is $\left(1+\varepsilon^{-1}\right)$--quasiharmonic in the sense of \cite{LW16-harmonic}. Indeed, if $w\in W^{1,2}(\Omega', X)$ satisfies $\trace(w) = \trace(v|_{\Omega'})$ then the map $w'$ which agrees with $w$ on $\Omega'$ and with $v$ on $\Omega\setminus \Omega'$ belongs to $W^{1,2}(\Omega, X)$ and satisfies $\trace(w') = \trace(u)$ by \cite[Theorem 1.12.3]{KS93}. Since $v$ minimizes $A_\varepsilon$ we obtain $$\Area(v|_{\Omega'}) + \varepsilon\cdot E_+^2(v|_{\Omega'}) \leq \Area(w) + \varepsilon\cdot E_+^2(w) \leq (1+\varepsilon)\cdot E_+^2(w)$$ and this implies the claim.

Finally, since $v$ is quasiharmonic and has a continuous trace, it follows from \cite[Theorem 1.3]{LW16-harmonic} that $v$ has a continuous representative which continuously extends to the boundary. This representative, which we denote again by $v$, satisfies the properties in the statement of the proposition.
\end{proof}

\bp\label{prop:Newton-Sobolev-rep}
 Let $\Omega$ be a Lipschitz domain in the interior of $M$. A measurable and essentially separably valued map $u\colon\Omega\to X$ belongs to $W^{1,2}(\Omega, X)$ if and only if there exist a map $v\colon\overline{\Omega}\to X$ and a Borel function $\rho\colon\overline{\Omega}\to[0,\infty]$ in $L^2(\Omega)$ such that $v=u$ almost everywhere and 
 \begin{equation}\label{eq:upper-grad-ineq}
  d(v(\gamma(a)),v(\gamma(b)))\leq \int_a^b \rho(\gamma(t))|\gamma'(t)|{\rm d} t
 \end{equation}
  for every Lipschitz curve $\gamma\colon[a,b]\to \overline{\Omega}$. In this case, we have
  \begin{equation}\label{eq:energy}
  E_+^2(u,g) = \inf\{\|\rho\|_{L^2(\Omega, g)}^2:\ \rho \text{ satisfies \eqref{eq:upper-grad-ineq}}\}
  \end{equation}
  and $v(z) = \trace(u)(z)$ for $\hm^1_g$--almost every $z\in\partial \Omega$.
\ep

The map $v$ in the claim is called a \emph{Newton--Sobolev representative} of $u$, and $\rho$ an \emph{upper gradient} of $v$. Inequality \eqref{eq:upper-grad-ineq} is known as the upper gradient inequality.

\begin{proof}
	The existence of $v$ and $\rho$ as in the claim imply that  $u\in W^{1,2}(\Omega,X)$, see \cite[Chapter 7]{HKST15}. For the opposite implication, by possibly doubling $M$ we may assume $M$ has no boundary. Since $\partial \Omega$ is Lip\-schitz, there exists a Lipschitz domain $\widehat \Omega\subset M$ containing $\overline\Omega$ and a map $\hat u\in W^{1,2}(\widehat \Omega,X)$ with $\hat u|_\Omega=u$, see the proof of \cite[Lemma 3.4]{LW15-Plateau}. There exists $v\colon \widehat \Omega\to X$ and $\rho\colon \widehat \Omega\to [0,\infty]$ satisfying \eqref{eq:upper-grad-ineq} for all Lipschitz curves $\gamma\colon [a,b]\to \widehat\Omega$, cf. \cite[Theorems 7.1.20 and 7.4.5]{HKST15}. The maps $v|_{\overline\Omega}$ and $\rho|_{\overline\Omega}$ satisfy the claim.
	
	The equality \eqref{eq:energy} follows e.g. from \cite[Theorem 7.1.20 and Lemma 6.2.2]{HKST15}. Let $\psi\colon (0,1)\times [0,1)\to M$ be as in the definition of the trace, so that $\trace(u)(\psi(s,0))=\lim_{t\to 0}u\circ\psi(s,t)$ for a.e. $s\in (0,1)$. A Fubini-type argument shows that  $$v(\psi(s,0))=\lim_{t\to 0}u\circ\psi(s,t)$$  for a.e. $s\in (0,1)$. This completes the proof.
\end{proof}

We illustrate the use of Newton-Sobolev representatives in the next lemma. Recall that a metric space is said to be $C$--quasiconvex if any two points can be joined by a Lipschitz curve of length at most $C$ times their distance.
\bl\label{lem:ug}
Let $h\colon K\to M$ be a Lipschitz map from a cell-complex $K$, and $A\subset K^1$ a $C$--quasiconvex subset of the 1-skeleton. If $v\colon M\to X$ is a Newton-Sobolev representative, and $\rho\in L^2(M)$ an upper gradient of $u$ with $ L:=\int_A\rho^2\circ h{\rm d} \hm^1<\infty,$ then $v\circ h|_A$ is $\frac 12$--H\"older continuous with constant $(CL)^{\frac 12}{\rm Lip}(h)$.
\el
\begin{proof}
	For $x,y\in A$, let $\gamma\colon [0,\ell(\gamma)]\to A$ be a simple unit speed curve joining $x$ and $y$ with $\ell(\gamma)\le Cd(x,y)$. By \eqref{eq:upper-grad-ineq} we have
	\begin{align*}
	d(v\circ h(x),v\circ h(y))&\le \int_0^{\ell(\gamma)}\rho\circ h(\gamma(t))|(h\circ\gamma)'(t)|{\rm d} t
	\le {\rm Lip}(h)\int_0^{\ell(\gamma)}\rho\circ h(\gamma(t)){\rm d}t\\
	&\le\  {\rm Lip}(h)\ell(\gamma)^{\frac 12}L^{\frac 12}\le (CL)^{\frac 12}{\rm Lip}(h)d(x,y)^{\frac 12}.
	\end{align*}
\end{proof}

\section{Admissible deformations on a surface}\label{sec:admissible-deformations}

The notion of admissible deformation on a surface given below, in the spirit of \cite{HL05}, will be used to define $1$--homotopy classes relative to given Jordan curves. We remark that the deformations in \cite{HL05, White88} keep the boundary fixed and are thus not suitable for studying the Plateau-Douglas problem. The deformations in \cite{White86, White88, HL03} for closed surfaces also do not adapt to our purposes.

\bd\label{def:admdef}
 An admissible deformation on a surface $M$ is a smooth map $\Phi\colon M\times \R^m\to M$, for some $m\in\N$, such that $\Phi_\xi:= \Phi(\cdot, \xi)$ is a diffeomorphism for every $\xi\in\R^m$ and $\Phi_0=\operatorname{id}_M$, and such that the derivative of $\Phi^p:= \Phi(p, \cdot)$ at the origin satisfies 
  \begin{equation*}
 D\Phi^p(0)(\R^m) = \left\{\begin{array}{l@{\;\;\text{if}\;}l} 
  T_pM & p\in {\rm int}(M)\\
  T_p(\partial M) &p\in\partial M.
 \end{array}\right.
 \end{equation*}
\ed

If $\Phi\colon M\times\R^m\to M$ is an admissible deformation on $M$ and $\varphi\colon M\to M$ is a diffeomorphism then $\Phi'(p,\xi):= \varphi(\Phi(\varphi^{-1}(p),\xi))$ also defines an admissible deformation.

\bp\label{prop:existence-admissible-deformations}
 There exist admissible deformations on every surface.
\ep

\begin{proof}
Let $\eta_1, \eta_2\colon [0,\infty)\to [0,\infty)$ be smooth functions such that $\eta_1(0)=0$, $\eta_1'(0)>0$, $\eta_2(0)>0$ and $\eta_1(t)=\eta_2(t) = 0$ for all $t\geq 1$. We use $\eta_1,\eta_2$ to define smooth vector fields $X_1, X_2$ on $M$ as follows. Each boundary component of $M$ has a neighborhood which is diffeomorphic to $S^1\times[0,2)$. On such a boundary component we define $X_1$ and $X_2$ by $X_1(z,t) = \eta_1(t) \frac{\partial}{\partial t}$ and $X_2(z,t) = \eta_2(t)\frac{\partial }{\partial z}$, written in coordinates $(z,t)\in S^1\times[0,2)$. Now, extend $X_1,X_2$ to all of $M$ by setting them to be zero outside these neighborhoods. It is easy to see that there exist smooth vector fields $X_3, \dots, X_m$ on $M$, for some $m$, with support in the interior of $M$ such that the vectors $X_1(p), \dots, X_m(p)$ span $T_pM$ for every $p$ in the interior of $M$. For every $k=1,\dots, m$, the flow $\varphi_{X_k, t}$ along $X_k$ is defined for all times $t\in\R$. Now the map $\Phi\colon M\times\R^m\to M$ given by $$\Phi(p,\xi):= \varphi_{X_1, \xi_1}\circ \varphi_{X_2, \xi_2}\circ\dots\circ\varphi_{X_m, \xi_m}(p)$$ defines an admissible deformation on $M$.
\end{proof}

Let $M$ be a surface, which we equip with a Riemannian metric $g$. Let $\Phi\colon M\times \R^m\to M$ be an admissible deformation on $M$ and let $h\colon K\to M$ be a triangulation of $M$. For $\xi\in\R^m$ let $h_\xi\colon K\to M$ be the triangulation given by $h_\xi:= \Phi_\xi\circ h$. The following variant of \cite[Lemma 5]{HL05} plays a key role throughout the article (see also \cite[Lemma 3.3]{HL03} for closed manifolds).

\bl\label{lem:crucial-ineq-deformation}
 There exist an open ball $B_{\Phi, h}\subset \R^m$ centered at the origin and $C>0$ such that for every Borel function $\rho\colon M\to [0,\infty]$ we have
   \begin{equation}\label{eq:first-ineq-crucial-lemma}
  \int_{B_{\Phi, h}}\left(\int_{K^0\cap \partial K} \rho\circ h_\xi(z)\,d\hm^0(z)\right)\,d\xi \leq C\int_{\partial M}\rho\,d\hm^1_g
 \end{equation}
and, for every $l\in\{0,1,2\}$,
 \begin{equation}\label{eq:second-ineq-crucial-lemma}
  \int_{B_{\Phi, h}}\left(\int_{K^l\setminus \partial K} \rho\circ h_\xi(z)\,d\hm^l(z)\right)\,d\xi \leq C\int_M\rho\,d\hm^2_g.
 \end{equation}
\el

\begin{proof}
 We only prove \eqref{eq:second-ineq-crucial-lemma} and leave the similar proof of \eqref{eq:first-ineq-crucial-lemma} to the reader. Let $\Delta$ be a closed cell of some dimension $l$ in $K$ and suppose $\Delta$ is not contained in $\partial K$. Define a map $H\colon \Delta\times\R^m\to M$ by $H(z,\xi):= \Phi(h(z), \xi)$. The properties of $\Phi$ and $h$ imply that $$D H(z, 0)(\R^l\times\R^m) = T_{h(z)}M$$ for every $z\in \Delta$ and therefore there exist $\varepsilon, c>0$ such that the jacobian of the differential of $H$ satisfies
 \begin{equation}\label{eq:lower-bound-jac}
 \jac(D H(z,\xi)) \geq c
 \end{equation} for every $(z,\xi)\in \Delta\times \overline{B}(0, 2\varepsilon)$. 
Since $H|_{\Delta\times B(0, 2\varepsilon)}$ is $C^1$ up to the boundary, we may extend $H$ to a map $ H\colon \tilde\Delta\times B(0, 2\varepsilon)\to \tilde M$ satisfying \eqref{eq:lower-bound-jac} for some open manifolds $\tilde\Delta\subset \R^l$ and $\tilde M$ containing $\Delta$ and $M$, respectively, by possibly making $c$ smaller.

We now claim that there exists $L\geq 0$ such that 
 \begin{equation}\label{eq:upper-bound-levelsets}
 \hm^{l+m-2}(H^{-1}(x) \cap \Delta\times B(0,\varepsilon))\leq L
\end{equation}
for every $x\in M$. In order to prove this, fix $(z,\xi)\in \Delta\times \overline{B}(0,\varepsilon)$. Let $F\colon \tilde\Delta\times B(0,2\varepsilon) \to \R^{l+m-2}$ be a $C^1$ map such that $F(z,\xi)=0$ and such that the map $$\tilde{H}\colon \tilde \Delta\times B(0,2\varepsilon)\to \tilde M\times\R^{l+m-2}$$ given by $\tilde{H} = (H, F)$ satisfies
$$D\tilde{H}(z,\xi) (\R^l\times\R^m) = T_{H(z,\xi)}\tilde M\times \R^{l+m-2}.$$
There exist $\delta>0$ and open neighborhoods $U\subset \tilde \Delta\times B(0,2\varepsilon)$ of $(z,\xi)$ and $V\subset \tilde M$ of $H(z,\xi)$ such that the restriction of $\tilde{H}$ to $U$ is a biLipschitz homeomorphism with image $V\times B(0,\delta)$. Let $G$ be the inverse of $\tilde{H}|_U$, so that $$H^{-1}(x)\cap U = G(\{x\}\times B(0,\delta))$$ for every $x\in V$. It follows that there exists $L'$ such that $$\hm^{l+m-2}(H^{-1}(x)\cap U) \leq L'$$ for every $x\in V$. Since $\Delta\times\overline{B}(0,\varepsilon)$ is compact we can cover it by finitely many such open sets $U$ and the claim follows for a suitable number $L$.
 
Finally, let $\rho\colon M\to[0,\infty]$ be a Borel function. From the co-area formula and the inequalities \eqref{eq:lower-bound-jac} and \eqref{eq:upper-bound-levelsets} we conclude
\begin{equation*}
 \begin{split}
  \int_{B(0,\varepsilon)}\int_\Delta \rho\circ H(z, \xi)&\,d\hm^l(z)\,d\xi\\
   &\leq c^{-1} \int_{B(0,\varepsilon)}\int_\Delta \rho\circ H(z, \xi)\, \jac(D H(z,\xi))\,d\hm^l(z)\,d\xi\\
  & = c^{-1}\int_M\rho(x) \cdot \hm^{l+m-2}(H^{-1}(x))\,d\hm^2_g(x)\\
  &\leq \frac{L}{c} \int_M\rho(x)\,d\hm^2_g(x).
 \end{split}
\end{equation*}
This proves \eqref{eq:second-ineq-crucial-lemma} with $B_{\Phi, h}:= B(0,\varepsilon)$ and $C=\frac{L}{c}$.
\end{proof}

Lemma~\ref{lem:crucial-ineq-deformation} has the following immediate corollary.

\bc\label{cor:aeagree}
If $N\subset M$ and $E\subset \partial M$ satisfy $\hm^2_g(N)=\hm^1_g(E)=0$ then, for almost every $\xi\in B_{\Phi,h}$, we have that $h_\xi(x)\notin N$ for $\hm^1$--a.e. $x\in K^1\setminus \partial K$ and $h_\xi(x)\notin E$ for every $x\in K^0\cap \partial K$.
\ec

In the next statement, we denote by $u\circ h_\xi|_{K^1}$ the  map which agrees with $u\circ h_\xi$ on $K^1\setminus\partial K$ and with $\trace(u)\circ h_\xi$ on $\partial K$.

\bp\label{prop:restriction-1-skeleton}
Let $X$ be a complete metric space and $u\in W^{1,2}(M,X)$. Then $u\circ h_\xi|_{K^1\setminus\partial K}$ is essentially continuous for a.e. $\xi\in B_{\Phi,h}$. If $u$ has continuous trace then $u\circ h_\xi|_{K^1}$ is essentially continuous for a.e. $\xi\in B_{\Phi,h}$, and extends continuously to $K$ in case $X$ is proper and admits a local quadratic isoperimetric inequality.
\ep

\begin{proof}
Let $v\colon M\to X$ be a Newton-Sobolev representative of $u$ with upper gradient $\rho\in L^2(M)$ (cf. Proposition~\ref{prop:Newton-Sobolev-rep}), and $A:=\overline{K^1\setminus \partial K}$. Since
$$
\int_A\rho^2\circ h_\xi{\rm d}\hm^1<\infty
$$
for a.e. $\xi\in B_{\Phi,h}$ by \eqref{eq:second-ineq-crucial-lemma}, Lemma~\ref{lem:ug} implies that $v\circ h_\xi|_A$ is continuous for a.e. $\xi\in B_{\Phi,h}$. The first claim now follows from Corollary~\ref{cor:aeagree} applied to the null-set $\{u\ne v\}$.

Suppose $\trace(u)$ has a continuous representative $\eta$. The set $\{v|_{\partial M}\ne \eta\}\subset \partial M$ has null $\hm^1_g$--measure by Proposition~\ref{prop:Newton-Sobolev-rep}, in particular $v\circ h_\xi|_{\partial K}=\eta\circ h_\xi|_{\partial K}$ $\hm^1$--a.e., for a.e. $\xi$. The argument above together with Corollary~\ref{cor:aeagree} applied to $\{u\ne v\}$ and $\{v|_{\partial M}\ne \eta\}$ implies that, for a.e. $\xi\in B_{\Phi,h}$, the map
\begin{align*}
w_\xi(x):= \left\{
\begin{array}{ll}
v\circ h_\xi(x), & x\in K^1\setminus\partial K\\
\eta\circ h_\xi(x), & x\in \partial K
\end{array}\right.
\end{align*}
is a continuous representative of $u\circ h_\xi|_{K^1}$. If $X$ is proper and admits a local quadratic isoperimetric inequality and $\Delta$ is a $2$--cell of $K$ then $\trace(u\circ h_\xi|_\Delta)=w_\xi|_{\partial\Delta}$ for a.e. $\xi\in B_{\Phi,h}$ by Proposition~\ref{prop:Newton-Sobolev-rep}. Applying Proposition~\ref{prop:good-cont-filling} on each $2$--cell and gluing these together yields the desired continuous extension.
\end{proof}

Now, suppose that $M$ has $k\geq 1$ boundary components and let $\Gamma$ be the disjoint union of $k$ rectifiable Jordan curves in a proper metric space $X$ admitting a local quadratic isoperimetric inequality. Recall the definition of homotopy relative to $\Gamma$ from the introduction.

\bt\label{thm:homotopic-1-skeleton}
Let $u\in\Lambda(M,\Gamma, X)$. Then there exists a negligible set $N\subset B_{\Phi, h}$ such that the continuous representatives of $u\circ h_\xi|_{K^1}$ and $u\circ h_\zeta|_{K^1}$ are homotopic relative to $\Gamma$ for all $\xi, \zeta\in B_{\Phi, h}\setminus N$.
\et

\begin{proof}
Denote by $\eta$ the continuous representative of $\trace(u)$. Let $v\colon M\to X$ be a Newton-Sobolev representative of $u$ with upper gradient $\rho\in L^2(M)$ as in Proposition~\ref{prop:Newton-Sobolev-rep}, and set $\bar v:=v$ on ${\rm int}(M)$ and $\bar v:=\eta$ on $\partial M$.	By the proof of Proposition~\ref{prop:restriction-1-skeleton}, there exists a null-set $N_0\subset B_{\Phi,h}$ such that $\bar v\circ h_\xi|_{K^1}$ is the continuous representative of $u\circ h_\xi|_{K^1}$ whenever $\xi\in B_{\Phi,h}\setminus N_0$.

We claim that there exists $\xi_0\in B_{\Phi,h}\setminus N_0$ such that the map $H_\xi\colon K^1\times [0,1]\to M$ given by $H_\xi(x,t):=\Phi(h(x),\xi_0+t(\xi-\xi_0)$ satisfies 
\begin{align}\label{eq:riesz}
\int_0^1\int_{K^l\setminus \partial K}\rho^2\circ H_\xi{\rm d}\hm^l{\rm d}t<\infty,\quad l=0,1,
\end{align}
for a.e. $\xi\in B_{\Phi,h}$.
Let us first finish the proof assuming \eqref{eq:riesz}. It is enough to show that there exists a null-set $N\subset B_{\Phi,h}$ containing $N_0$ such that $\bar v\circ h_{\xi_0}|_{K^1}\sim \bar v\circ h_{\xi}|_{K^1}$ rel $\Gamma$, whenever $\xi\in B_{\Phi,h}\setminus N$. Indeed, from this it follows that $\bar v\circ h_{\xi}|_{K^1}\sim \bar v\circ h_{\zeta}|_{K^1}$ rel $\Gamma$ for every $\xi,\zeta\in B_{\Phi,h}\setminus N$.

Note that $H_\xi(\cdot,0)=\bar v\circ h_{\xi_0}|_{K^1}$, $H_\xi(\cdot,1)=\bar v\circ h_\xi|_{K^1}$ and that $\bar v\circ H_\xi|_{\partial K\times [0,1]}$ is continuous with $\bar v\circ H_\xi|_{\partial K\times \{t\}}\in [\Gamma]$ for every $\xi\in B_{\Phi,h}$ and $t\in [0,1]$. Fix a $1$--cell $e$ of $K$ not contained in $\partial K$ and let $A:=e\times [0,1]$. We show that $\bar v\circ H_\xi|_{\partial A}$ is continuous and the trace of a Sobolev map, for a.e $\xi\in B_{\Phi,h}\setminus N_0$. By Proposition~\ref{prop:good-cont-filling} this implies that $\bar v\circ H_\xi|_{\partial A}$ has a continuous extension to $A$, and choosing a continuous extension for each $A$ we obtain the desired homotopy relative to $\Gamma$ between $\bar v\circ h_{\xi_0}|_{K^1}$ and $\bar v\circ h_\xi|_{K^1}$, for a.e. $\xi\in B_{\Phi,h}\setminus N_0$.

Since ${\rm Lip}(H_\xi)\cdot\rho\circ H_\xi|_A$ is an upper gradient of $v\circ H_\xi|_A$, it follows from \eqref{eq:riesz} and Lemma~\ref{lem:ug} that $v\circ H_\xi|_A\in W^{1,2}(A,X)$ and $\trace(v\circ H_\xi|_A)=v\circ H_\xi|_{\partial A}$ for a.e. $\xi\in B_{\Phi,h}\setminus N_0$. For a.e. $\xi\in B_{\Phi,h}\setminus N_0$, we have that $\bar v\circ H_\xi(z_0,\cdot)=v\circ H_\xi(z_0,\cdot)$ is H\"older continuous for $z_0\in K^0\setminus \partial K$ by \eqref{eq:riesz} and Lemma~\ref{lem:ug}, and
$$
\bar v\circ H_\xi(z_0,t)= v\circ H_\xi(z_0,t)\quad \textrm{ a.e. }t\in [0,1]
$$
for $z_0\in K^0\cap \partial K$, by  Corollary~\ref{cor:aeagree} and a Fubini-type argument. Thus $\bar v\circ H_\xi|_{\partial A}$ is the continuous representative of $v\circ H_\xi|_{\partial A}$ for a.e. $\xi\in B_{\Phi,h}\setminus N_0$. This completes the proof that $\bar v\circ H_\xi|_{\partial A}$ is continuous and the trace of a Sobolev function, for a.e. $\xi\in B_{\Phi,h}\setminus N_0$.

It remains to show \eqref{eq:riesz}. Define
$$ f(\xi):=\chi_{B_{\Phi,h}}(\xi)\left(\int_{K^0\setminus\partial K}\rho^2\circ h_\xi{\rm d}\hm^0+\int_{K^1\setminus\partial K}\rho^2\circ h_\xi{\rm d}\hm^1\right),\quad \xi\in \R^m.$$
Then $f\in L^1(\R^m)$ by \eqref{eq:second-ineq-crucial-lemma} and thus there exists $\xi_0\in B_{\Phi,h}\setminus N_0$ such that the Riesz potential $R_1f(\xi_0):=\int_{\R^m}\frac{f(\xi_0+\xi)}{|\xi|^{m-1}}{\rm d}\xi$ is finite (cf. \cite[Theorem 3.22]{hei01}). Integrating in spherical  coordinates we have $\int_{S^{m-1}}\int_0^\infty f(\xi_0+tw){\rm d}t{\rm d}w= R_1f(\xi_0)<\infty.$
Since 
$$\int_0^1\int_{K^l\setminus \partial K}\rho^2\circ H_\xi{\rm d}\hm^l{\rm d}t\le \int_0^\infty f(\xi_0+t(\xi-\xi_0))){\rm d}t= |\xi-\xi_0|\int_0^\infty f\left(\xi_0+s\frac{\xi-\xi_0}{|\xi-\xi_0|}\right) {\rm d}s$$ 
for $l=0,1$ and $\xi\in B_{\Phi,h}\setminus\{\xi_0\}$, \eqref{eq:riesz} follows.
\end{proof}

We end this section with the following lemma which will be used in the proofs of the theorems in the next section.

\bl\label{lem:conv-restr-1-skeleton}
 Let $u\in W^{1,2}(M, X)$ and let $(u_n)\subset W^{1,2}(M,X)$ be an energy bounded sequence converging to $u$ in $L^2(M,X)$. Then for almost every $\xi\in B_{\Phi, h}$ there exists a subsequence $(u_{n_j})$ such that the continuous representative of $u_{n_j}\circ h_\xi|_{K^1\setminus \partial K}$ converges uniformly to the continuous representative of $u\circ h_\xi|_{K^1\setminus\partial K}$ as $j\to\infty$.
\el

\begin{proof}
	By passing to a subsequence we may assume that $u_n\to u$ almost everywhere in $M$. For each $n\in\N$, let $v_n\colon M\to X$ be a Newton-Sobolev representative of $u_n$ with upper gradient $\rho_n\in L^2(M)$ satisfying
		$$ \|\rho_n\|_{L^2(M,g)}^2\le 2E_+^2(u_n,g), $$
	cf. Proposition~\ref{prop:Newton-Sobolev-rep}. 
By the proof of Proposition~\ref{prop:restriction-1-skeleton} and Corollary~\ref{cor:aeagree}, there exists a negligible set $N_0\subset B_{\Phi, h}$ such that for every $z\in B_{\Phi,h}\setminus N_0$ the map $v_n\circ h_\xi|_{K^1\setminus \partial K}$ is the continuous representative of $u_n\circ h_\xi|_{K^1\setminus \partial K}$ for every $n\in\N$ and 
\begin{equation}\label{eq:vn-tou}
 v_n\circ h_\xi|_{K^1\setminus \partial K} \rightarrow u\circ h_\xi|_{K^1\setminus \partial K}
\end{equation}
 $\hm^1$--a.e. with $n\to\infty$.
Set $A:=\overline{K^1\setminus\partial K}$. Fatou's lemma and \eqref{eq:second-ineq-crucial-lemma} imply that
	\begin{align*}
	\int_{B_{\Phi,h}}\left(\liminf_{n\to\infty}\int_A\rho_n^2\circ h_\xi{\rm d}\hm^1\right){\rm d}\xi&\le 	\liminf_{n\to\infty}\int_{B_{\Phi,h}}\int_A\rho_n^2\circ h_\xi{\rm d}\hm^1{\rm d}\xi\\
	&\le C\liminf_{n\to\infty}\int_M\rho_n^2{\rm d}\hm^2_g<\infty.
	\end{align*}
Therefore, for almost every $\xi\in B_{\Phi, h}\setminus N_0$, we have $$\liminf_{n\to\infty}\int_A\rho_n^2\circ h_\xi{\rm d}\hm^1<\infty.$$ By Lemma~\ref{lem:ug},  Arzela-Ascoli's Theorem and  \eqref{eq:vn-tou}, for such $\xi$ there exists a subsequence $(v_{n_j}\circ h_\xi|_A)_{j\in\N}$ which is uniformly $\frac 12$--H\"older continuous and converges uniformly to the continuous representative of $u\circ h_\xi|_{K^1\setminus \partial K}$ as $j\to\infty$.
\end{proof}

\section{The relative $1$--homotopy class of Sobolev maps}\label{sec:1-homot}

Throughout this section, let $X$ be a proper geodesic metric space admitting a local quadratic isoperimetric inequality. Let $\Gamma\subset X$ be the disjoint union of $k\geq 1$ rectifiable Jordan curves, and let $M$ be a surface with $k$ boundary components. We fix a Riemannian metric $g$ on $M$.

Let $\Phi\colon M\times\R^m\to M$ be an admissible deformation on $M$. Theorem~\ref{thm:homotopic-1-skeleton} shows that for every $u\in\Lambda(M,\Gamma, X)$ and every triangulation $h\colon K\to M$ of $M$ we have $$[u\circ h_\xi|_{K^1}]_\Gamma = [u\circ h_\zeta|_{K^1}]_\Gamma$$ for almost all $\xi,\zeta\in B_{\Phi, h}$. We denote the common relative homotopy class by $u_{\#, 1}[h]$. The following theorem shows that $u_{\#,1}[h]$ is  independent of the choice of deformation $\Phi$ and that inducing the same relative homotopy class is independent of the triangulation $h$.

\bt\label{thm:rel-hom-indep-wiggling}
 Let $X$, $\Gamma$, $M$, $\Phi$ be as above. Let $u\in \Lambda(M,\Gamma, X)$ and let $h\colon K\to M$ be a triangulation of $M$. The relative homotopy class $u_{\#,1}[h]$ does not depend on the choice of admissible deformation $\Phi$. Moreover, if $v\in\Lambda(M,\Gamma, X)$ is such that $v_{\#,1}[h] = u_{\#,1}[h]$ then we have $v_{\#,1}[\tilde{h}] = u_{\#,1}[\tilde{h}]$ for any triangulation $\tilde{h}\colon \tilde{K}\to M$.
\et

We will  need the following two lemmas in the proof.

\bl\label{lem:approx-cont-area-energy-bounded}
 Let $u\in W^{1,2}(M, X)$ have continuous trace. Then for all $\varepsilon, \delta>0$ there exists a continuous map $\hat{u}\colon M\to X$ in $W^{1,2}(M,X)$ with $\hat{u}|_{\partial M} = \trace(u)$, $d_{L^2}(u, \hat{u})<\varepsilon$, and 
 \begin{equation}\label{eq:bound-area-energy-cont}
  \Area(\hat{u}) \leq \Area(u) + \delta \cdot E_+^2(u,g),\quad E_+^2(\hat{u}, g) \leq \left(1+\delta^{-1}\right) \cdot E_+^2(u,g).
 \end{equation}
\el

\begin{proof}
Let $u$ be as in the statement of the lemma and let $\varepsilon, \delta>0$. Fix an admissible deformation $\Phi$ on $M$ and let $\varepsilon'>0$ be sufficiently small, to be determined later. Choose a triangulation $h\colon K\to M$ of $M$ in such a way that for every $\xi\in B_{\Phi, h}$ we have $\hm^2_g(h_\xi(\Delta))<\varepsilon'$ for every $2$--cell $\Delta\subset K$. 
  
It follows from (the proof of) Proposition~\ref{prop:restriction-1-skeleton} that, for almost every $\xi\in B_{\Phi, h}$, the map $u\circ h_\xi|_{K^1}$ is essentially continuous and its restriction to the boundary of each open $2$--cell $\Delta\subset K$ coincides with the trace of the Sobolev map $u\circ h_\xi|_{\Delta}$. Fix such $\xi$ and abbreviate $H:= h_\xi$. It thus follows that if $\Delta$ is an open $2$--cell then the map $u|_{H(\partial \Delta)}$ is essentially continuous and the trace of the Sobolev map $u|_{H(\Delta)}$. By Proposition~\ref{prop:good-cont-filling} there thus exists a continuous map $u_\Delta\colon H(\overline{\Delta})\to X$ which extends the continuous representative of $u|_{H(\partial \Delta)}$, belongs to $W^{1,2}(H(\Delta), X)$ and satisfies
$$\Area(u_\Delta) \leq \Area(u|_{H(\Delta)}) + \delta\cdot E_+^2(u|_{H(\Delta)}, g)$$ 
as well as 
$$E_+^2(u_\Delta, g) \leq \left(1+\delta^{-1}\right)\cdot E_+^2(u|_{H(\Delta)}, g).$$

It follows from the Sobolev-Poincar\'e inequality (see \cite[Section 2]{Heb99} for closed manifolds), from  \cite[Corollary 1.6.3]{KS93} and H\"older's inequality that 
\begin{equation*}
 \begin{split}
   \int_{H(\Delta)} d^2(u_\Delta(z), u(z))\,d\hm^2_g(z)&\leq C\cdot \hm_g^2(H(\Delta))\cdot \left[E_+^2(u_\Delta, g) + E_+^2(u|_{H(\Delta)}, g)\right]\\
   &\leq C\varepsilon' \left(2+\delta^{-1}\right) \cdot E_+^2(u|_{H(\Delta)}, g)
 \end{split}
\end{equation*}
for some constant $C$ depending on $(M,g)$.

Finally, let $\hat{u}\colon M\to X$ be the continuous map obtained by gluing the maps $u_\Delta$ along their boundaries. Then $\hat u\in W^{1,2}(M, X)$ by \cite[Theorem 1.12.3]{KS93} and, taking the sum over all $\Delta$ in the three inequalities above, we obtain the inequalities in \eqref{eq:bound-area-energy-cont} as well as $$\int_{M} d^2(\hat{u}(z), u(z))\,d\hm^2_g(z)\leq C\varepsilon' \left(2+\delta^{-1}\right) \cdot E_+^2(u,g).$$ Upon choosing $\varepsilon'>0$ sufficiently small, this yields $d_{L^2}(\hat{u}, u)<\varepsilon$.
\end{proof}

\bl\label{lem:close-1-skeleton-implies-homotopic-rel}
 Let $X$, $\Gamma$, $M$ be as above. Then there exists $\delta>0$ with the following property. Let $h\colon K\to M$ be a triangulation and let $\varrho,\varrho'\colon K^1\to X$ be continuous such that $\varrho|_{\partial K},\varrho'|_{\partial K}\in [\Gamma]$ are homotopic via a family of maps in $[\Gamma]$. If $$\sup_{z\in K^1\setminus \partial K} d(\varrho(z), \varrho'(z))< \delta$$ and if for every component $C$ of $\partial M$ for which the Jordan curve $\varrho(C)$ is not contractible in $X$ we have 
 \begin{equation}\label{eq:bound-on-boundary}
  \sup_{z\in C} d(\varrho(z), \varrho'(z))< \delta
 \end{equation}
   then $\varrho$ and $\varrho'$ are homotopic relative to $\Gamma$.
\el

The condition \eqref{eq:bound-on-boundary} cannot be omitted, as easy examples show. The lemma will also be used in the proof of Theorem~\ref{thm:stability-1-homotopic}, where it will be essential that we do not impose any condition akin to \eqref{eq:bound-on-boundary} for the components $C$ of $\partial K$ which are mapped to contractible Jordan curves.

\begin{proof}
 Since $X$ is proper, geodesic and admits a local quadratic isoperimetric inequality it follows from \cite[Theorem 5.2]{LWY20}, \cite[Proposition 2.2]{LWY20}, and from the proof of \cite[Proposition 6.2]{LWY20} that there exists $r_0>0$ such that every closed curve in $X$ of diameter at most $4r_0$ is contractible. Recall that $\Gamma=\Gamma_1\cup\dots\cup \Gamma_k$ is the disjoint union of rectifiable Jordan curves. We may assume that $3r_0\leq\diam(\Gamma_i)$ for every $i$. There exists $0<\delta<r_0/3$ such that whenever $x,y\in \Gamma$ satisfy $d(x,y)\leq 9\delta$ then they belong to the same Jordan curve $\Gamma_i$ and one of the two segments of $\Gamma_i$ joining $x$ and $y$ has diameter at most $r_0$. 
 
Let $\varrho, \varrho'\colon K^1\to X$ be as in the statement of the lemma with this specific choice of $\delta$.  
After possibly adding vertices to $K^1\setminus \partial K$ we may further assume that the image under $\varrho$ and $\varrho'$ of any edge in the closure of $K^1\setminus \partial K$ has diameter at most $\delta$. 
We now construct a homotopy $H\colon K^1\times[0,1]\to X$ relative to $\Gamma$ between $\varrho$ and $\varrho'$. Let $H(\cdot, 0) = \varrho$ and $H(\cdot,1)=\varrho'$. For each $z_0\in K^0\setminus \partial K$ let $H(z_0,\cdot)$ be a (constant speed) geodesic from $\varrho(z_0)$ and $\varrho'(z_0)$. 
For each component $C$ of $\partial K$ and each $z\in K^0\cap C$, let $H(z,\cdot)$ be a weakly monotone parametrization of one of the segments in $\varrho(C)$ joining $\varrho(z)$ to $\varrho'(z)$ in such a way that, for every edge $e\subset C$, the map $H|_{\partial (e\times [0,1])}$ is contractible in $\varrho(C)$. By \eqref{eq:bound-on-boundary}, in the case that $\varrho(C)$ is not contractible, we may choose $H(z,\cdot)$ to have diameter at most $r_0$ for every $z\in C\cap K^0$.

For every edge $e\subset \partial K$ the map $H|_{\partial (e\times[0,1])}$ admits a continuous extension with image in $\Gamma$ such that for each $t\in[0,1]$ the map $H(\cdot, t)$ is weakly monotone. Moreover, for every edge $e\subset K^1$ not intersecting $\partial K$ the curve $H|_{\partial(e\times[0,1])}$ has diameter at most $4\delta$ and thus admits a continuous extension to $e\times [0,1]$. Finally, let $e\subset K^1$ be an edge which intersects (but is not contained in) some component $C$ of $\partial K$. Notice that the image of $H|_{\partial (e\times[0,1])}$ is contained in the $3\delta$--neighbourhood of $\varrho(C)$. Thus, if $\varrho(C)$ is contractible then $H|_{\partial (e\times[0,1])}$ admits a continuous extension to $e\times[0,1]$. If $\varrho(C)$ is not contractible then, by construction, the image of $H|_{\partial (e\times[0,1])}$ has diameter at most $4r_0$ and hence admits again a continuous extension to $e\times[0,1]$.
\end{proof}

\begin{proof}[Proof of Theorem~\ref{thm:rel-hom-indep-wiggling}]
 Let $u\in \Lambda(M,\Gamma, X)$ and let $h\colon K\to M$ be a triangulation.  We wish to show that the relative homotopy class, which we denote by $u_{\#, 1}[h,\Phi]$ for the moment, is independent of the choice of admissible deformation $\Phi$. Let $(u_n)$ be a sequence of continuous maps $u_n\colon M\to X$ converging in $L^2(M,X)$ to $u$, with $u_n|_{\partial M} = \trace(u)$ and $u_n\in W^{1,2}(M,X)$ for every $n$, and such that the energy of $u_n$ is bounded independently of $n$. Such a sequence exists by Lemma~\ref{lem:approx-cont-area-energy-bounded} and we call it a good approximating sequence for $u$. 
 
 We first claim that there exists a subsequence $(n_j)$ such that $u_{\#, 1}[h,\Phi] = [u_{n_j}\circ h|_{K^1}]_\Gamma$ for all $j\geq 1$. Indeed, by Proposition~\ref{prop:restriction-1-skeleton} and Theorem~\ref{thm:homotopic-1-skeleton} there exists a negligible subset $N\subset B_{\Phi, h}$ such that for $\xi,\xi'\in B_{\Phi, h}\setminus N$ the maps $u\circ h_\xi|_{K^1}$ and $u\circ h_{\xi'}|_{K^1}$ are essentially continuous and their continuous representatives are homotopic relative to $\Gamma$. Since $u_n|_{\partial M} = \trace(u)$ it follows with Lemma~\ref{lem:conv-restr-1-skeleton} that for almost every $\xi_0\in B_{\Phi, h}\setminus N$ there is a subsequence $(n_j)$ such that the maps $u_{n_j}\circ h_{\xi_0}|_{K^1}$ converge uniformly to the continuous representative of $u\circ h_{\xi_0}|_{K^1}$ as $j\to \infty$. Fix such $\xi_0$ and such a subsequence $(n_j)$. Lemma~\ref{lem:close-1-skeleton-implies-homotopic-rel} thus implies that there exists $j_0$ such that $u_{n_j}\circ h_{\xi_0}|_{K^1}$ is homotopic relative to $\Gamma$ to the continuous representative of $u\circ h_{\xi_0}|_{K^1}$ for every $j\geq j_0$. Since $u_{n_j}$ is continuous the maps $u_{n_j}\circ h_{\xi_0}|_{K^1}$ and $u_{n_j}\circ h|_{K^1}$ are homotopic relative to $\Gamma$. It thus follows that for all $j\geq j_0$ the continuous representative of $u\circ h_\xi|_{K^1}$ is homotopic relative to $\Gamma$ to $u_{n_j}\circ h|_{K^1}$ for every $\xi\in B_{\Phi, h}\setminus N$. Upon reindexing the subsequence we may assume that $j_0=1$. This proves the claim.
  
It easily follows from the claim that $u_{\#,1}[h,\Phi]$ is independent of $\Phi$. Indeed, let $\tilde{\Phi}$ be another admissible deformation on $M$. On the one hand, the claim shows that there exists a subsequence $(n_j)$ such that $$u_{\#, 1}[h,\Phi] = [u_{n_j}\circ h|_{K^1}]_\Gamma$$ for all $j\geq 1$. Applying the claim again with $\Phi$ replaced by $\tilde{\Phi}$ and with $(u_n)$ replaced by $(u_{n_j})$ we see that there is a further subsequence $(n_{j_l})$ such that $$u_{\#, 1}[h,\tilde{\Phi}] = [u_{n_{j_l}}\circ h|_{K^1}]_\Gamma$$ for all $l\geq 1$. From this it follows that $u_{\#,1}[h,\Phi] = u_{\#,1}[h,\tilde{\Phi}]$, which proves the first statement of the theorem.

The second statement of the theorem also follows from the claim. Indeed, let $v\in\Lambda(M,\Gamma, X)$ be such $v_{\#, 1}[h] = u_{\#,1}[h]$ and let $(v_n)$ be a good approximating sequence for $v$. The claim shows that we can find a subsequence $(n_j)$ such that $$[u_{n_j}\circ h|_{K^1}]_\Gamma = u_{\#,1}[h] = v_{\#,1}[h] = [v_{n_j}\circ h|_{K^1}]_\Gamma$$ for all $j\geq 1$. Let $\tilde{h}\colon \tilde{K}\to M$ be another triangulation. Since $u_{n_j}$ and $v_{n_j}$ are continuous it is easy to see that $$[u_{n_j}\circ\tilde{h}|_{\tilde{K}^1}]_\Gamma = [v_{n_j}\circ \tilde{h}|_{\tilde{K}^1}]_\Gamma$$ for all $j\geq 1$, compare with \cite[Lemma 2.1]{HL03}. The claim now implies that $u_{\#,1}[\tilde{h}] = v_{\#,1}[\tilde{h}]$, which proves the second statement of the theorem.
\end{proof}

\bp
 Let $\varphi\colon M\to X$ be a continuous map such that $\varphi|_{\partial M}\in[\Gamma]$ and let $u\in\Lambda(M, \Gamma, X)$. Then $$u_{\#,1}[h] = [\varphi\circ h|_{K^1}]_\Gamma$$ holds for one triangulation $h\colon K \to M$ if and only if it holds for every triangulation.
\ep

\begin{proof}
 Let $h\colon K\to M$ be a triangulation of $M$ such that $$u_{\#,1}[h] = [\varphi\circ h|_{K^1}]_\Gamma$$
 and let $(u_n)$ be a good approximating sequence for $u$ as in the first paragraph of the proof of Theorem~\ref{thm:rel-hom-indep-wiggling}. By the claim in the second paragraph of that proof, there exists a subsequence $(n_j)$ such that $u_{\#,1}[h] = [u_{n_j}\circ h|_{K^1}]_\Gamma$ for all $j\geq 1$ and hence $$[u_{n_j}\circ h|_{K^1}]\Gamma = [\varphi\circ h|_{K^1}]_\Gamma$$ for all $j\geq 1$. Let $\tilde{h}\colon\tilde{K}\to M$ be another triangulation of $M$. Since $u_{n_j}$ and $\varphi$ are continuous $$[u_{n_j}\circ \tilde{h}|_{\tilde{K}^1}]_\Gamma = [\varphi\circ\tilde{h}|_{\tilde{K}^1}]_\Gamma$$ for all $j\geq 1$, compare with \cite[Lemma 2.1]{HL03}. After possibly passing to a further subsequence we have $u_{\#,1}[\tilde{h}] = [u_{n_j}\circ \tilde{h}|_{\tilde{K}^1}]_\Gamma$ for all $j\geq 1$ and hence $$u_{\#,1}[\tilde{h}]=[\varphi\circ\tilde{h}|_{\tilde{K}^1}]_\Gamma.$$ This concludes the proof.
\end{proof}

\bd
 Two maps $u,v\in\Lambda(M,\Gamma, X)$ are said to be $1$--homotopic relative to $\Gamma$, denoted $u\sim_1 v$ rel $\Gamma$, if for some and thus every triangulation $h$ of $M$ we have $u_{\#,1}[h] = v_{\#,1}[h]$. If $u\in\Lambda(M,\Gamma, X)$ and $\varphi\colon M\to X$ is continuous with $\varphi|_{\partial M}\in[\Gamma]$ then $u$ and $\varphi$ are said to be $1$--homotopic relative to $\Gamma$, denoted $u\sim_1\varphi$ rel $\Gamma$, if for some and thus every triangulation $h\colon K\to M$ we have $u_{\#,1}[h] = [\varphi\circ h|_{K^1}]_\Gamma$.
\ed

If $u,v\in \Lambda(M,\Gamma,X)$, $u\sim_1 v$ rel $\Gamma$ and $\psi\colon M\to M$ is a diffeomorphism then $u\circ \psi\sim_1 v\circ\psi$ rel $\Gamma$, see the remark after Definition~\ref{def:admdef}.

\bt\label{thm:stability-1-homotopic}
 Let $X$, $\Gamma$, $M$ be as above. Then for every $L>0$ there exists $\varepsilon>0$ such that if $u,v\in\Lambda(M,\Gamma, X)$ induce the same orientation on $\Gamma$ and satisfy $$\max\left\{E_+^2(u,g), E_+^2(v,g)\right\}\leq L\quad\textrm{and}\quad d_{L^2}(u,v)\leq \varepsilon,$$ then $u$ and $v$ are $1$--homotopic relative to $\Gamma$. 
\et

Notice that the theorem does not imply the stability of $1$--homotopy classes relative to $\Gamma$, since the $L^2$--limit of a sequence in $\Lambda(M,\Gamma,X)$ with uniformly bounded energy need not belong to $\Lambda(M,\Gamma,X)$. An analog of Theorem~\ref{thm:stability-1-homotopic} holds for closed surfaces (where $\Gamma=\varnothing$ and $\Lambda(M,\Gamma,X)=W^{1,2}(M,X)$) and in this case implies the stability of $1$--homotopy classes in the presence of a local quadratic isoperimetric inequality. Example~\ref{ex:stabilityfail} below shows that the local quadratic isoperimetric inequality is crucial for this.

\begin{proof}
 We argue by contradiction and assume the statement is not true. Then there exist energy bounded sequences $(u_n),(v_n)\subset\Lambda(M,\Gamma, X)$ such that for every $n\in\N$ we have $d_{L^2}(u_n,v_n)\leq \frac{1}{n}$, that $u_n$ and $v_n$ induce the same orientation on $\Gamma$ but $u_n$ is not $1$--homotopic to $v_n$ relative to $\Gamma$. After possibly passing to a subsequence, we may assume by the Rellich-Kondrachov compactness theorem \cite[Theorem 1.13]{KS93} and by \cite[Lemma 2.4]{FW-Plateau-Douglas} that there exists $u\in W^{1,2}(M, X)$ such that the sequences $(u_n)$ and $(v_n)$ both converge to $u$ in $L^2(M, X)$. 
 
 Fix an admissible deformation $\Phi$ on $M$ and a triangulation $h\colon K\to M$. By Proposition~\ref{prop:restriction-1-skeleton} and Theorem~\ref{thm:homotopic-1-skeleton} there exists a negligible set $N\subset B_{\Phi, h}$ such that for all $\xi,\zeta\in B_{\Phi, h}\setminus N$ and all $n\in\N$ we have that $u_n\circ h_\xi|_{K^1}$ and $u_n\circ h_\zeta|_{K^1}$ are essentially continuous and their continuous representatives are homotopic relative to $\Gamma$ and that the same is true when $u_n$ is replaced by $v_n$ and $u$. It moreover follows from Lemma~\ref{lem:conv-restr-1-skeleton} that for almost every $\xi_0\in B_{\Phi, h}\setminus N$ there exists a subsequence $(n_j)$ such that the continuous representatives of $u_{n_j}\circ h_{\xi_0}|_{K^1\setminus \partial K}$ and of $v_{n_j}\circ h_{\xi_0}|_{K^1\setminus \partial K} $ both converge uniformly to the continuous representative of $u\circ h_{\xi_0}|_{K^1\setminus \partial K}$. Fix such $\xi_0$ and denote by $\varrho_j$ and $\varrho'_j$ the continuous representatives of  $u_{n_j}\circ h_{\xi_0}|_{K^1}$ and of $v_{n_j}\circ h_{\xi_0}|_{K^1}$, respectively. Denote by $C_m$ and $\Gamma_m$, $m=1,\dots, k$, the components of $\partial K$ and $\Gamma$, respectively. Notice that the sequences $(\varrho_j|_{\partial K})_j$ and $(\varrho'_j|_{\partial K})_j$ both converge in $L^2(\partial K, X)$ to $u\circ h_{\xi_0}|_{\partial K}$ by \cite[Theorem 1.12.2]{KS93}. Thus, after possibly relabelling the components, we may assume that $$\varrho_j(C_m) = \Gamma_m = \varrho'_j(C_m)$$ for all sufficiently large $j$ and every $m=1,\dots, k$. Since $\varrho_j|_{\partial K}$ and $\varrho'_j|_{\partial K}$ induce the same orientation on $\Gamma$ it follows, in particular, that $\varrho_j|_{\partial K}$ and $\varrho'_j|_{\partial K}$ are homotopic via a family of maps in $[\Gamma]$ for every sufficiently large $j$. Let $m$ be such that $\Gamma_m$ is not contractible. Then it follows from the remark after Proposition~\ref{prop:good-cont-filling} together with the proof of \cite[Proposition 5.1]{FW-Plateau-Douglas} that the families $(\varrho_j|_{C_m})$ and $(\varrho'_j|_{C_m})$ are both equi-continuous. Hence, after possibly passing to a subsequence, we may assume that both sequences converge uniformly to the continuous representative of $u\circ h_{\xi_0}|_{C_m}$. It thus follows that for every sufficiently large $j$ the maps $\varrho_j$ and $\varrho'_j$ satisfy the hypotheses of Lemma~\ref{lem:close-1-skeleton-implies-homotopic-rel}. In particular, it follows that there exists $j_0$ such that $\varrho_j$ and  $\varrho_j'$ are homotopic relative to $\Gamma$ for every $j\geq j_0$. Hence, for every $\xi\in B_{\Phi, h}\setminus N$ and every $j\geq j_0$ we have that the continuous representatives of  $u_{n_j}\circ h_{\xi}|_{K^1}$ and $v_{n_j}\circ h_{\xi}|_{K^1}$ are homotopic relative to $\Gamma$. This shows that $u_{n_j}$ and $v_{n_j}$ are $1$--homotopic relative to $\Gamma$, which is a contradiction, concluding the proof.
\end{proof}

\begin{proof}[Proof of Theorem~\ref{thm:intro-properties-1-hom-class-Sobolev}]
	Statements (ii) and (iii) follow from Theorems~\ref{thm:rel-hom-indep-wiggling} and \ref{thm:stability-1-homotopic}. As for statement (i), suppose $u$ has a continuous representative $\bar{u}\colon M\to X$. We have  $u\circ h_\xi|_{K^1} = \bar{u}\circ h_\xi|_{K^1}$ $\hm^1$--a.e., for almost every $\xi$ by Corollary~\ref{cor:aeagree} and hence $$[u\circ h_\xi|_{K^1}]_\Gamma =[\bar u\circ h_\xi|_{K^1}]_\Gamma= [\bar{u}\circ h|_{K^1}]_\Gamma$$ for almost every $\xi$. This proves statement (i).
\end{proof}

\begin{example}\label{ex:stabilityfail}
Consider the surface of revolution $C\subset \R^3$ of the graph of $$f\colon (0,1]\to [1/3,1], \quad  f(x)= (2+\sin(1/x))/3.$$ The compact set  $C\cup \{0\}\times \overline{\mathbb D}\subset \R^3$ equipped with the subspace metric is not geodesic, but by adding a countable number of suitable line segments parallel to the $x$--axis, connecting points on $C$ to $\{0\}\times \overline{\mathbb D}$, we obtain a compact subset of $\R^3$ bi-Lipschitz equivalent to a geodesic space $Y$. It is not difficult to see that $Y$, and thus $X:=S^1\times Y$, fails to admit a local quadratic isoperimetric inequality.
Let $x_n\to 0$ be the sequence of local minima of $f$, and $h_n\colon S^1\to Y$ the constant speed parametrizations corresponding to the circles $\{x_n\}\times \R^2\cap C$. The maps $$u_n\colon S^1\times S^1\to X,\quad (z,z')\mapsto (z,h_n(z'))$$ are bi-Lipschitz for each $n$, and converge  uniformly to the map $u(z,z')=(z,h(z'))$, where $h\colon S^1\to Y$ is the constant speed parametrization of the circle corresponding to $\{(0,z'/3):\ z'\in S^1\}\subset Y$.
However, one can check that the maps $h_n$ are all non-contractible and pairwise $1$--homotopic, while $h$ is contractible. It follows that $u$ cannot lie in the common homotopy class of the maps $u_n$.
\end{example}

The example above can be modified so that the maps $u_n$ form an area minimizing sequence in their common $1$--homotopy class. Considering the set $C\cup \{0\}\times\overline{\mathbb D}$ with the metric inherited from $\R^3$ in the example above, we obtain a \emph{non-geodesic} space with a local quadratic isoperimetric inequality where the stability of $1$--homotopy classes of maps from closed surfaces fails.

\section{The homotopic Douglas condition and its consequences}\label{sec:homot-Douglas-cond}

Let $X$ be a proper geodesic metric space admitting a local quadratic isoperimetric inequality, and let $\Gamma\subset X$ be the disjoint union of $k\geq 1$ rectifiable Jordan curves. Let $M$ be a connected surface with $k$ boundary components, and let $\varphi\colon M\to X$ be a continuous map such that $\varphi|_{\partial M}\in[\Gamma]$.

\bp\label{prop:induced-homom-Douglas}
 If the induced homomorphism $\varphi_*\colon \pi_1(M)\to \pi_1(X)$ on fundamental groups is injective then $\varphi$ satisfies the homotopic Douglas condition \eqref{eq:Douglas-condition}.
\ep

\begin{proof}
 We first claim that $a(M,\varphi, X)<\infty$. Let $l_0>0$ be as in the definition of the local quadratic isoperimetric inequality. Since $\Gamma$ is a finite union of rectifiable Jordan curves there exists $0<r_0<l_0/3$ such that every subcurve of $\Gamma$ of diameter at most $r_0$ has length at most $l_0/3$. Moreover we may choose $r_0$ small enough so that all closed loops of diameter $\le 2r_0$ are contractible, cf. the proof of Lemma~\ref{lem:close-1-skeleton-implies-homotopic-rel}. Now, fix a triangulation of $M$ all of whose $2$--cells are triangles. We identify the $1$--skeleton of the triangulation with a subset of $M$ and denote it by $M^1$. Choosing the triangulation sufficiently fine we may assume that for each $1$--cell $e\subset M^1$ we have $\diam(\varphi(e))< r_0$. Let $u\colon M^1\to X$ be the continuous map which agrees with $\varphi$ on the $0$--skeleton $M^0$ and such that for each $1$--cell $e\subset M^1$ the following holds: if $e$ is contained in $\partial M$ then $u|_e$ is the constant speed parametrization of $\varphi(e)$; if $e$ is not contained in $\partial M$ then $u|_e$ is a geodesic. It follows that for every $2$--cell $\Delta\subset M$ the curve $u|_{\partial \Delta}$ is Lipschitz and has length at most $l_0$ and thus has a continuous Sobolev extension to $\Delta$ (which we denote $u|_{\Delta}$) by the local quadratic isoperimetric inequality and Lemma~\ref{lem:approx-cont-area-energy-bounded}. Also note that $u|_e$ is end-point homotopic to $\varphi|_e$ by the choice of $r_0$. The continuous map $\bar u\colon M\to X$ obtained by gluing all the $u|_{\Delta}$ together is a Sobolev map and satisfies $\bar u|_{M^1} \sim \varphi|_{M^1}$ rel $\Gamma$. It thus follows that $\bar u\sim_1\varphi$ relative to $\Gamma$. The map $\bar u$ has finite area and thus we obtain $a(M,\varphi, X)<\infty$, as claimed.
 
Since the induced homomorphism $\varphi_*\colon \pi_1(M)\to \pi_1(X)$ on fundamental groups is injective it follows that if $\alpha$ is a simple closed non-contractible curve in the interior of $M$ then $\varphi\circ \alpha$ is not contractible. Consequently, there are no primary reductions $(M^*, \varphi^*)$ of $(M,\varphi)$ and hence $a^*(M,\varphi,X)=\infty$ by definition. Since $a(M,\varphi, X)<\infty$ this shows that $\varphi$ satisfies the homotopic Douglas condition.
\end{proof}

\bp\label{prop:equi-cont}
 Let $g$ be a Riemannian metric on $M$. Then for every $\eta>0$ and $L>0$ the family $$\{\trace(u): \text{$u\in\Lambda(M,\Gamma, X)$, $u\sim_1 \varphi$ rel $\Gamma$, $E_+^2(u,g)\leq L$, $\Area(u)\leq a^*(M,\varphi, X) - \eta$}\}$$ is equi-continuous.
\ep

A corresponding result without fixing relative $1$--homotopy classes is contained in \cite[Proposition 5.1]{FW-Plateau-Douglas}. In order to control the relative $1$--homotopy class of the maps that we construct in the proof, we will use the following technical lemma. 

In the next statement, $\alpha$ is a smooth closed simple non-contractible curve in the interior of $M$ and let $M^*$ be the smooth surface obtained from $M$ by cutting $M$ along $\alpha$ and gluing smooth discs to the two newly created boundary components.

\bl\label{lem:1-hom-reduction}
 Let $A\subset M$ be a biLipschitz cylinder such that $A\cap \partial M$ is connected and one boundary component of $A$ coincides with $\alpha$. Suppose there is $v\in \Lambda(M^*, \Gamma, X)$ inducing the same orientation on $\Gamma$ as $u$ and satisfying $v|_{M\setminus A} = u$. Then $\varphi\circ \alpha$ is contractible and $v$ is $1$--homotopic to $\varphi^*$ relative to $\Gamma$, whenever $\varphi^*\colon M^*\to X$ is continuous and coincides with $\varphi$ on $M\setminus \alpha$.
\el

\begin{proof}
 Let $A'\subset M$ be a biLipschitz cylinder with piecewise smooth boundary components and such that $A'$ contains a small neighborhood of $A$ in $M$. The boundary component $\alpha'$ of $A'$ which is homotopic to $\alpha$ outside $A$ is contained in the interior of $M$. Let $\beta'$ be the other boundary component of $A'$. If $\gamma:= A\cap \partial M$ is not empty then $\beta'$ contains $\gamma$. 
 
Let $h\colon K\to M$ be a triangulation of $M$ such that $h(K^1)$ contains $\alpha'$ and $\beta'$. Let $K'$ be the sub-complex of $K$ obtained by removing the interior of cells that get mapped to the interior of $A'$. Let $K^*$ be the complex obtained from $K'$ by adding two cells, each glued along the preimage of $\alpha'$ and $\beta'$, respectively, and extend $h|_{K'}$ to a triangulation $h^*\colon K^*\to M^*$ of $M^*$. Let $C\subset K'^1$ be the preimage of $\beta'\cap\partial M$ under $h$.

Let $U\subset M$ be a small neighborhood of $\alpha$ whose closure is contained in the interior of $A'$. Using vector fields as in the proof of Proposition~\ref{prop:existence-admissible-deformations} it is not difficult to construct admissible deformations $\Phi\colon M\times\R^m\to M$ on $M$ and $\Phi^*\colon M^*\times \R^m\to M^*$ on $M^*$ which agree on $(M\setminus U)\times B(0,\varepsilon)$ for some sufficiently small $\varepsilon>0$. On $K'^1\setminus C$ the maps $h_\xi=\Phi_\xi\circ h$ and $h^*_\xi = \Phi^*_\xi\circ h^*$ agree for sufficiently small $\xi$ and stay outside $A$, so we have $$v\circ h^*_\xi|_{K'^1\setminus C} = u\circ h_\xi|_{K'^1\setminus C}$$ for a.e. small $\xi$. 
Since $u$ and $v$ induce the same orientation on $\Gamma$ it follows that $v\circ h^*_\xi|_{K'^1}$ is homotopic to $u\circ h_\xi|_{K'^1}$ relative to $\Gamma$ for almost every sufficiently small $\xi$.

Now, $u\circ h_\xi|_{K'^1}$ is homotopic to $\varphi\circ h|_{K'^1}$ relative to $\Gamma$ for almost every $\xi$ sufficiently small. Let $\Omega\subset M^*$ be the Lipschitz Jordan domain bounded by $\alpha'$. Since $v\circ\Phi^*_\xi|_{\partial \Omega}$ is the trace of the Sobolev disc $v\circ\Phi^*_\xi|_{\Omega}$ for almost every small $\xi$ it follows from Proposition~\ref{prop:good-cont-filling} that the continuous representative of $v\circ\Phi^*_\xi\circ \alpha'$ is contractible and hence $\varphi\circ\alpha'$ and therefore $\varphi\circ\alpha$ are contractible. Let $\varphi^*\colon M^*\to X$ be a continuous extension of $\varphi|_{M\setminus \alpha}$ to $M^*$. Since $$\varphi^*\circ h^*|_{K'^1} = \varphi\circ h|_{K'^1}$$ and the $1$--skeletons of $K^*$ and $K'$ agree it follows that $v\circ h^*_\xi|_{K^{*1}}$ is homotopic to $\varphi^*\circ h^*|_{K^{*1}}$ relative to $\Gamma$ for almost every $\xi$ sufficiently small. This shows that $v$ is $1$--homotopic to $\varphi^*$ relative to $\Gamma$.  
\end{proof}

The proof of Proposition~\ref{prop:equi-cont} is almost the same as that of \cite[Proposition 5.1]{FW-Plateau-Douglas}, so we only give a rough sketch.

\begin{proof}[Proof of Proposition~\ref{prop:equi-cont}]
Denote by $\mathscr A$ the family of maps $u\in \Lambda(M,\Gamma,X)$ such that $u\sim_1 \varphi$ rel $\Gamma$,  $E^2_+(u,g)\le L$ and $\Area(u)\le a^*(M,\varphi,X)-\eta$. Suppose the claim is not true. Then there exists $\varepsilon_0>0$ and, for each $\delta>0$, a map $u\in \mathscr A$ such that the image of some boundary arc with length $\le \delta$ has length $\ge \varepsilon_0$. By considering a conformal chart containing the short boundary arc and using the Courant-Lebesgue lemma \cite[Lemma 7.3]{LW15-Plateau} we see that there exists an arc $\beta\colon I\to M$ connecting two boundary points on either side (and outside) of the short boundary arc, for which $u\circ \beta\in W^{1,2}(I,X)$ agrees with the continuous representative of $\trace(u)$ at the end-points, and $\length(u\circ\beta)\le \pi[E^2_+(u,g)/\log(1/\delta)]^{1/2}$.

Since $\Gamma$ consists of rectifiable Jordan curves, there exists $\delta'>0$ so that any points on $\Gamma$ with distance at most $\delta'$ belong to the same component and the shorter of the arcs joining them  has length $<\min\{\varepsilon_0,\eta'\}$, where $0<\eta'<l_0/2$ is such that $C(2\eta')^2<\eta/2$. Here $C$ and $l_0$ are the constants in the local quadratic isoperimetric inequality of $X$. Thus, by choosing $\delta>0$ small enough, it follows that $\length(u\circ\beta)<\eta'$ and moreover the image $\Gamma^+$ of the longer boundary arc $\gamma^+$ joining the endpoints of $\beta$  has length $<\eta'$.

Let $\alpha\subset \inter M$ be a smooth Jordan curve bounding  an annulus $A\subset M$ together with the curve $\alpha':=\gamma^+\cup\beta$ such that $u\circ\alpha\in W^{1,2}(S^1,X)$. In the surface $M^*$ obtained by cutting $M$ along $\alpha$ and gluing discs to the newly created boundary curves, $\alpha'$ bounds a Lipschitz Jordan domain $\Omega$. If $\Gamma_0$ is the concatenation of $u\circ\beta$ and $\Gamma^+=\trace(u)\circ\gamma^+$, then $\length(\Gamma_0)<2\eta'$ and, by \cite[Lemma 4.8]{LW-intrinsic}, $\Gamma_0$ is the trace of a Sobolev map $w_\Omega\in W^{1,2}(\Omega,X)$ with $\Area(w_\Omega)< C(2\eta')^2<\eta/2$.

We define $v$ as $w_{\Omega}$ and $u|_{M\setminus A}$ on the respective sets. To define $v$ on the remaining smooth disc $\Omega'\subset M^*$, map $A$ diffeomorphically to an annulus $A'\subset \Omega'$ identifying $\alpha$ with $\partial\Omega'$, and $\alpha'$ with a Jordan curve (compactly contained in $\Omega'$) that bounds a copy $\Omega''$ of $\Omega$, and set $v|_{\Omega''}=w_\Omega$ and $v|_{A'}=u|_A$ (after the diffeomorphic identifications). The gluing theorem \cite[Theorem 1.12.3]{KS93} implies that $v\in W^{1,2}(M^*,X)$ and by construction $v\in \Lambda(M^*,\Gamma,X)$ with $v$ and $u$ inducing the same orientation on $\Gamma$. Lemma~\ref{lem:1-hom-reduction} implies that $v$ is $1$--homotopic to $\varphi^*$ rel $\Gamma$ for any primary reduction $(M^*,\varphi^*)$ of $(M,\varphi)$. Now the estimate
\[
\Area(v)=\Area(u|_{M\setminus A})+2\Area(w_\Omega)+\Area(u|_A)<\Area(u)+\eta
\]
yields a contradiction with the fact that $u\in\mathscr A$, completing the proof.
\end{proof}

In the next proposition, we assume that the Euler characteristic $\chi(M)$ of $M$ is strictly negative so that $M$ admits a hyperbolic metric, that is, a Riemannian metric on $M$ of constant curvature $-1$ and such that $\partial M$ is geodesic. 

\bp\label{prop:lower-bound-rel-systole}
For every $\eta>0$ and $L>0$ there exists $\varepsilon>0$ with the following property. If $u\in\Lambda(M,\Gamma, X)$ is $1$--homotopic to $\varphi$ relative to $\Gamma$ and such that $$\Area(u)\leq a^*(M,\varphi, X)-\eta,$$ and if $g$ is a hyperbolic metric on $M$ satisfying $E_+^2(u,g)\leq L$ then the relative systole of $(M,g)$ satisfies $\sys_{\rm rel}(M,g)\geq \varepsilon$.
\ep

The relative systole $\sys_{\rm rel}(M,g)$ of $(M,g)$ is the minimal length of curves $\beta$ in $M$ of the following form. Either $\beta$ is closed and not contractible in $M$ via a family of closed curves, or the endpoints of $\beta$ lie on the boundary of $M$ and $\beta$ is not contractible via a family of curves with endpoints on $\partial M$. 
The proof of the proposition is almost the same as that of \cite[Proposition 6.1]{FW-Plateau-Douglas} and we only sketch it. Lemma~\ref{lem:1-hom-reduction} will be used again to control the relative $1$--homotopy type of the primary reductions appearing in the proof.

\begin{proof}
Let $\beta_0$ be the geodesic realizing the systole $\lambda:=\sys_{\rm rel}(M,g)$. We may use a collar neighbourhood to find a 'parallel' Jordan curve $\beta\colon I\to M$ for which $u\circ\beta\in W^{1,2}(I,X)$ and $\ell(u\circ\beta)\le 2[\lambda E^2_+(u,g)]^{1/2}$, see \cite[Lemma 6.2]{FW-Plateau-Douglas}. If $\beta$ connects two boundary points, then $I$ is a closed interval and the proof is analogous to that of Proposition~\ref{prop:equi-cont}. Namely, using the notation from the proof of Proposition~\ref{prop:equi-cont} and supposing the relative systole $\lambda$ is small enough, we may assume the boundary points are on the same boundary component and the image $\Gamma^+$ of one boundary arc $\gamma^+$ connecting them has small length, so that the concatenation $\Gamma_0$ of $u\circ\beta$ and $\Gamma^+$ satisfies $\ell(\Gamma_0)<2\eta'$.

We let $\alpha\subset \inter M$ be a closed Jordan curve bounding a (closed) annulus $A$ with $\alpha':=\gamma^+\cup\beta$ such that $u\circ\alpha\in W^{1,2}(S^1,X)$. In the surface $M^*$ obtained from $M$ by cutting along $\alpha$, $\alpha'$ bounds a Jordan domain $\Omega$ containing $A$ and we let $w_\Omega\in W^{1,2}(\Omega,X)$ satisfy $\trace(w_\Omega)=\Gamma_0$ and $\Area(w_\Omega)< C(2\eta')^2<\eta/2$. Defining $v\in \Lambda(M^*,\Gamma,X)$ as in the proof of Proposition~\ref{prop:equi-cont}, we reach the same contradiction with the fact that $\Area(u)\le a^*(M,\varphi,X)-\eta$.

If $\beta_0$ is a closed geodesic, we construct $M^*$ and $v$ essentially as in the proof of \cite[Proposition 6.1]{FW-Plateau-Douglas} (keeping any components without boundary, and defining $v$ on them analogously). We omit the details.
\end{proof}

\section{Solution of the homotopic Plateau-Douglas problem}\label{sec:sol}

Let $X$ be a proper geodesic metric space admitting a local quadratic isoperimetric inequality and let $\Gamma\subset X$ be the union of $k\geq 1$ rectifiable Jordan curves. Let $M$ be a connected surface with $k$ boundary components and let $\varphi\colon M\to X$ be a continuous map such that $\varphi|_{\partial M}\in [\Gamma]$. 

\bp\label{prop:tech-min-seq-1-hom}
Suppose $\chi(M)<0$. Let $(u_n)\subset\Lambda(M, \Gamma, X)$ be a sequence such that each $u_n$ is $1$--homotopic to $\varphi$ relative to $\Gamma$ and $$\sup_n \Area(u_n)<a^*(M,\varphi, X).$$ Let $(g_n)$ be a sequence of hyperbolic metrics on $M$. Then there exist $u\in\Lambda(M,\Gamma, X)$ which is $1$--homotopic to $\varphi$ relative to $\Gamma$ and a hyperbolic metric $g$ on $M$ such that $$\Area(u)\leq \limsup_{n\to\infty} \Area(u_n)\quad\text{ and }\quad E_+^2(u,g)\leq \limsup_{n\to\infty} E_+^2(u_n, g_n).$$
\ep

\begin{proof}
 Let $(u_n)$ and $(g_n)$ be as in the statement of the proposition. By \cite[Theorem 1.2 and (5.2)]{FW-Morrey} there exist hyperbolic metrics $\tilde{g}_n$ such that $$E_+^2(u_n, \tilde{g}_n) \leq \frac{4}{\pi}\cdot \Area(u_n) +1.$$ After possibly replacing $g_n$ by $\tilde{g}_n$ and passing to a subsequence, we may therefore assume that the energies $E_+^2(u_n, g_n)$ are uniformly bounded and converge to a limit denoted by $m$.

By Proposition~\ref{prop:lower-bound-rel-systole}, the relative systoles of $(M, g_n)$ are uniformly bounded away from zero. Therefore, by the Mumford compactness theorem (see \cite[Theorem 3.3]{FW-Plateau-Douglas} and \cite[Theorem 4.4.1]{DHT10} for the fact that the diffeomorphisms may be chosen to be orientation preserving), there exist orientation preserving diffeomorphisms $\psi_n\colon M\to M$ and a hyperbolic metric $h$ on $M$ such that, after possibly passing to a subsequence, the Riemannian metrics $\psi_n^*g_n$ smoothly converge to $h$. For $n\in\N$ define a map by $v_n:= u_n\circ\psi_n$ and notice that $v_n\in\Lambda(M,\Gamma, X)$. Since $\psi_n$, when viewed as a map from $(M, h)$ to $(M, g_n)$, is $\lambda_n$-biLipschitz with $\lambda_n\to 1$ it follows that $$\lim_{n\to\infty} E_+^2(v_n, h) = m.$$  By \cite[Lemma 2.4]{FW-Plateau-Douglas} and the metric space valued Rellich-Kondrachov theorem (see \cite[Theorem 1.13]{KS93}) there exists $v\in W^{1,2}(M,X)$ such that a subsequence $(v_{n_j})$ converges in $L^2(M,X)$ to $v$. The lower semi-continuity of energy implies that $E_+^2(v, h)\leq m$. Since each $u_n$ is $1$--homotopic to $\varphi$ relative to $\Gamma$ and each $\psi_n$ is orientation preserving it follows that all the maps $v_n$ induce the same orientation on $\Gamma$. 
By Theorem~\ref{thm:stability-1-homotopic} there thus exists $j_0\in\N$ such that $v_{n_j}$ is $1$--homotopic to $v_{n_{j_0}}$ for every $j\geq j_0$. It follows that for $j\geq j_0$ the maps $w_j:= v_{n_j}\circ \psi_{n_{j_0}}^{-1}\in \Lambda(M,\Gamma, X)$ satisfy $$w_j\sim_1 v_{n_{j_0}}\circ \psi_{n_{j_0}}^{-1}=u_{n_{j_0}} \sim_1 \varphi\textrm{ rel }\Gamma.$$ The sequence $(w_j)$ converges in $L^2(M,X)$ to the map $u:= v\circ\psi_{n_{j_0}}^{-1}$ and $g:=(\psi_{n_{j_0}}^{-1})^*h$, we furthermore have $$E_+^2(u,g)\leq \lim_{j\to \infty} E_+^2(w_j, h_0)= \lim_{j\to\infty} E_+^2(v_j, h) = m.$$ Finally, Proposition~\ref{prop:equi-cont} implies that the family $\{\trace(w_j): j\in\N\}$ is equi-continuous and hence, after passing to a further subsequence, we may assume that the sequence $(\trace(w_j))$ converges uniformly to some continuous map $\gamma\colon \partial M\to X$. As the uniform limit of weakly monotone parametrizations of $\Gamma$, the map $\gamma$ is also a weakly monotone parametrization of $\Gamma$. Since $(\trace(w_j))$ converges in $L^2(\partial M, X)$ to $\trace(u)$ it follows that $\trace(u)=\gamma$ and hence $u\in \Lambda(M,\Gamma, X)$. Since $u$ and $w_j$ induce the same orientation on $\Gamma$ and since $w_j$ is $1$--homotopic to $\varphi$ relative to $\Gamma$ for every $j$ sufficiently large, it follows from Theorem~\ref{thm:stability-1-homotopic} that $u$ is $1$--homotopic to $\varphi$ relative to $\Gamma$ as well. The lower semi-continuity of area and invariance of area under diffeomorphisms imply that $$\Area(u) \leq \liminf_{j\to\infty}\Area(w_j) \leq \limsup_{n\to\infty} \Area(u_n).$$ This concludes the proof.
\end{proof}

\begin{proof}[Proof of Theorem~\ref{thm:Plateau-Douglas-homot-intro}]
Let $X$, $M$, $\Gamma$ be as in the statement of the theorem and let $\varphi\colon M\to X$ be a continuous map with $\varphi|_{\partial M}\in[\Gamma]$ satisfying the Douglas condition \eqref{eq:Douglas-condition}. 

We start by proving (i) in the case $\chi(M)<0$. The family $$\Lambda_{\rm min}:= \{u\in\Lambda(M,\Gamma, X): \text{$u \sim_1 \varphi$ relative to $\Gamma$ and $\Area(u) = a(M,\varphi, X)$}\}$$ is not empty. Indeed, this follows from Proposition~\ref{prop:tech-min-seq-1-hom}, applied to a sequence $(u_n)\subset\Lambda(M,\Gamma, X)$ and an arbitrary sequence of hyperbolic metrics such that $u_n$ is $1$--homotopic to $\varphi$ relative to $\Gamma$ for every $n$ and $$\Area(u_n) \to a(M, \varphi, X)$$ as $n$ tends to infinity.
Next, set $$m:= \inf\{E_+^2(u, g): \text{$u\in \Lambda_{\rm min}$, $g$ hyperbolic metric}\}$$ and choose sequences $(u_n)$ and $(g_n)$, where $u_n\in \Lambda_{\rm min}$ and where the $g_n$ are hyperbolic metrics on $M$, such that $$\lim_{n\to\infty} E_+^2(u_n, g_n)= m.$$ Applying Proposition~\ref{prop:tech-min-seq-1-hom} to these sequences we obtain a map $u\in \Lambda_{\rm min}$ and a hyperbolic metric $g$ on $M$ such that $E_+^2(u,g) = m$. It now follows from \cite[Corollary 1.3]{FW-Morrey} that $u$ is infinitesimally isotropic with respect to $g$. 

We are left with the case $\chi(M)\ge 0$. If $k=1$ and ${\rm genus}(M)=0$, the result follows from \cite[Theorem 1.2 and 1.4]{FW-Plateau-Douglas} since in this case any two maps inducing the same orientation on $\Gamma$ are $1$--homotopic.

In the remaining case $k=2$ and ${\rm genus}(M)=0$, one uses the Mumford compactness theorem for flat metrics normalized to have volume 1 (see \cite[Theorem 4.4.1]{DHT10} for the case of closed surfaces) and a flat collar lemma to prove an analog of Proposition~\ref{prop:lower-bound-rel-systole}. Replacing Proposition~\ref{prop:lower-bound-rel-systole} by this analog, the proof of Proposition~\ref{prop:tech-min-seq-1-hom} remains valid, and the argument above then works verbatim. See the proof of \cite[Theorem 1.2]{FW-Plateau-Douglas} for more discussion. This concludes the proof of statement (i).

To show (ii) let $u$ and $g$ be as in statement (i). Then $u$ is a local area minimizer and it follows from the proof of \cite[Theorem 1.4]{FW-Plateau-Douglas} that $u$ has a representative $\bar{u}$ which is locally H\"older continuous in the interior of $M$ and continuously extends to the boundary $\partial M$, thus proving statement (ii).

Statement (iii) is a direct consequence of the following lemma.
\end{proof}

\bl\label{lem:pi2}
 Let $X$ be a metric space, let $\Gamma\subset X$ be the the disjoint union of $k\geq 1$ Jordan curve, and let $M$ be a smooth compact surface with $k$ boundary components. If $X$ has trivial second homotopy group then two continuous maps $\varphi, \psi\colon M\to X$ with $\varphi|_{\partial M}, \psi|_{\partial M}\in [\Gamma]$ are $1$--homotopic relative to $\Gamma$ if and only if they are homotopic relative to $\Gamma$.
\el

We provide the easy proof for completeness, compare with \cite[Lemma 2.1]{Lem82}.

\begin{proof}
 Let $X$, $M$, $\Gamma$ be as in the statement of the lemma and let $\varphi, \psi\colon M\to X$ be continuous maps such that $\varphi|_{\partial M}, \psi|_{\partial M}\in [\Gamma]$. It is clear that if $\varphi$ and $\psi$ are homotopic relative to $\Gamma$ then they are, in particular, $1$--homotopic relative to $\Gamma$. In order to prove the opposite direction, suppose $\varphi$ and $\psi$ are $1$--homotopic relative to $\Gamma$ and let $F\colon K^1\times [0,1]\to X$ be a homotopy from $\varphi$ to $\psi$ such that $F(\cdot, t)\in[\Gamma]$ for all $t$. Let $G$ be the continuous map which coincides with $F$ on $K^1\times [0,1]$ and with $\varphi$ and $\psi$ on $K\times\{0\}$ and $K\times\{1\}$, respectively. For every $2$--cell $\Delta\subset K$ the restriction of $G$ to $\partial(\Delta\times[0,1])$ extends to a continuous map on $\Delta\times[0,1]$ since $X$ has trivial second homotopy group. The map $\bar{G}\colon K\times[0,1]\to X$ obtained in this way is a homotopy relative to $\Gamma$ between $\varphi$ and $\psi$.
\end{proof}

Observe that being $1$--homotopic is a more restrictive condition than inducing the same action on fundamental groups.

\begin{example}\label{ex:1-hom-stronger-action-fundgrp}
	Let $X=S^1\times S^1$ be the standard torus, $\Gamma=\{1\}\times S^1\cup\{e^{i\pi}\}\times S^1\subset X$, and $M=[0,1]\times S^1$. The maps $\varphi_\pm\in \Lambda(M,\Gamma,X)$ given by $\varphi_\pm(t,z)=(e^{\pm i\pi t},z)$ induce the same action $\pi_1(M)\to \pi_1(X)$ and agree on $\partial M$, but are not $1$--homotopic relative to $\Gamma$. Note that $\varphi_\pm$ are both conformal area minimizers in $\Lambda(M,\Gamma,X)$.
\end{example}

We finish the paper by discussing an analog of Theorem~\ref{thm:Plateau-Douglas-homot-intro} for closed surfaces, that is, $k=0$. In this case $\Gamma=\varnothing$ and consequently $\trace(u)\in [\Gamma]$ is a vacuous condition; in particular $\Lambda(M,\Gamma,X)=W^{1,2}(M,X)$. We say that two maps are $1$--homotopic if they are $1$--homotopic relative to $\Gamma=\varnothing$.

We assume throughout this discussion that {\bf $X$ is compact}, so that the Rellich-Kondrachov compactness theorem is applicable for any energy bounded sequence in $W^{1,2}(M,X)$. (The assumption $\trace(u)\in [\Gamma]$ prevents a sequence from escaping to infinity when $\Gamma\ne \varnothing$, and we prevent the same here by assuming  compactness.) Thus the results in Section~\ref{sec:1-homot} about $1$--homotopy remain valid with these interpretations. Note that, with the convention $\sys_{\rm rel}(M)=\sys(M)$, Proposition~\ref{prop:lower-bound-rel-systole} (and thus Proposition~\ref{prop:tech-min-seq-1-hom}) also remain valid with the same proofs.

The following theorem extends \cite[Theorem 4.4]{SU82} and \cite[Theorem 3.1]{SY79} to non-smooth target spaces. 

\bt\label{thm:area-min-hom-class-without-bdry}
Suppose $M$ is a closed surface, and $X$ a compact geodesic metric space admitting a local quadratic isoperimetric inequality. If a continuous map $\varphi\colon M\to X$ satisfies the homotopic Douglas condition, then there exist $u\in W^{1,2}(M,X)$ and a Riemannian metric $g$ on $M$ such that $u$ is $1$--homotopic to $\varphi$, $u$ is infinitesimally isotropic with respect to $g$, and $$\Area(u) = a(M, \varphi, X).$$ 
Furthermore, any such $u$ has a representative $\bar{u}$ which is H\"older continuous in $M$. If $X$ has trivial second homotopy group then $\bar{u}$ is homotopic to $\varphi$ relative to $\Gamma$.
\et

\begin{proof}
The proof Theorem~\ref{thm:Plateau-Douglas-homot-intro} (as well as that of Lemma~\ref{lem:pi2}) remains valid under the hypotheses of the claim (see the discussion above), except for the existence of $u$ and $g$ in the case $\chi(M)\ge 0$, i.e.~$M=S^2$ or $M=S^1\times S^1$. 

In the first case we may choose $u\equiv {\rm constant}$ and $g$ the standard metric on $S^2$, since $\varphi$ is $1$--homotopic to a constant map. In the second case $M=S^1\times S^1$ we use Mumford's compactness theorem for flat metrics with volume normalized to 1 to obtain analogs of Propositions~\ref{prop:lower-bound-rel-systole} and \ref{prop:tech-min-seq-1-hom} and proceed as in the proof of Theorem~\ref{thm:Plateau-Douglas-homot-intro}.
\end{proof}

\def\cprime{$'$} \def\cprime{$'$} \def\cprime{$'$}


\begin{thebibliography}{10}

\bibitem{Cou40}
R.~Courant.
\newblock The existence of minimal surfaces of given topological structure
  under prescribed boundary conditions.
\newblock {\em Acta Math.}, 72:51--98, 1940.

\bibitem{Cou37}
Richard Courant.
\newblock Plateau's problem and {D}irichlet's principle.
\newblock {\em Ann. of Math. (2)}, 38(3):679--724, 1937.

\bibitem{Cre-Plateau-sing}
Paul Creutz.
\newblock Plateau's problem for singular curves.
\newblock {\em Comm. Anal. Geom., to appear}.

\bibitem{Creutz-Fitzi}
Paul Creutz and Martin Fitzi.
\newblock The {P}lateau--{D}ouglas problem for singular configurations and in
  general metric spaces.
\newblock {\em preprint ArXiv:2008.08922}, 2020.

\bibitem{DHT10}
Ulrich Dierkes, Stefan Hildebrandt, and Anthony~J. Tromba.
\newblock {\em Global analysis of minimal surfaces}, volume 341 of {\em
  Grundlehren der Mathematischen Wissenschaften [Fundamental Principles of
  Mathematical Sciences]}.
\newblock Springer, Heidelberg, second edition, 2010.

\bibitem{Dou31}
Jesse Douglas.
\newblock Solution of the problem of {P}lateau.
\newblock {\em Trans. Amer. Math. Soc.}, 33(1):263--321, 1931.

\bibitem{Dou39}
Jesse Douglas.
\newblock Minimal surfaces of higher topological structure.
\newblock {\em Ann. of Math. (2)}, 40(1):205--298, 1939.

\bibitem{FW-Plateau-Douglas}
Martin Fitzi and Stefan Wenger.
\newblock Area minimizing surfaces of bounded genus in metric spaces.
\newblock {\em J. Reine Angew. Math.}
\newblock Published online April 16th 2020.

\bibitem{FW-Morrey}
Martin Fitzi and Stefan Wenger.
\newblock Morrey's {$\varepsilon$}-conformality lemma in metric spaces.
\newblock {\em Proc. Amer. Math. Soc.}, 148(10):4285--4298, 2020.

\bibitem{HL03}
Fengbo Hang and Fanghua Lin.
\newblock Topology of {S}obolev mappings. {II}.
\newblock {\em Acta Math.}, 191(1):55--107, 2003.

\bibitem{HL05}
Fengbo Hang and Fanghua Lin.
\newblock Topology of {S}obolev mappings. {IV}.
\newblock {\em Discrete Contin. Dyn. Syst.}, 13(5):1097--1124, 2005.

\bibitem{Heb99}
Emmanuel Hebey.
\newblock {\em Nonlinear analysis on manifolds: {S}obolev spaces and
  inequalities}, volume~5 of {\em Courant Lecture Notes in Mathematics}.
\newblock New York University, Courant Institute of Mathematical Sciences, New
  York; American Mathematical Society, Providence, RI, 1999.

\bibitem{hei01}
Juha Heinonen.
\newblock {\em Lectures on analysis on metric spaces}.
\newblock Universitext. Springer-Verlag, New York, 2001.

\bibitem{HKST15}
Juha Heinonen, Pekka Koskela, Nageswari Shanmugalingam, and Jeremy Tyson.
\newblock {\em Sobolev spaces on metric measure spaces}, volume~27 of {\em New
  Mathematical Monographs}.
\newblock Cambridge University Press, Cambridge, 2015.

\bibitem{Jost84}
J\"{u}rgen Jost.
\newblock {\em Harmonic mappings between {R}iemannian manifolds}, volume~4 of
  {\em Proceedings of the Centre for Mathematical Analysis, Australian National
  University}.
\newblock Australian National University, Centre for Mathematical Analysis,
  Canberra, 1984.

\bibitem{Jos85}
J{\"u}rgen Jost.
\newblock Conformal mappings and the {P}lateau-{D}ouglas problem in
  {R}iemannian manifolds.
\newblock {\em J. Reine Angew. Math.}, 359:37--54, 1985.

\bibitem{Jost94}
J\"{u}rgen Jost.
\newblock Equilibrium maps between metric spaces.
\newblock {\em Calc. Var. Partial Differential Equations}, 2(2):173--204, 1994.

\bibitem{Kar07}
M.~B. Karmanova.
\newblock Area and co-area formulas for mappings of the {S}obolev classes with
  values in a metric space.
\newblock {\em Sibirsk. Mat. Zh.}, 48(4):778--788, 2007.

\bibitem{Kir94}
Bernd Kirchheim.
\newblock Rectifiable metric spaces: local structure and regularity of the
  {H}ausdorff measure.
\newblock {\em Proc. Amer. Math. Soc.}, 121(1):113--123, 1994.

\bibitem{KS93}
Nicholas~J. Korevaar and Richard~M. Schoen.
\newblock Sobolev spaces and harmonic maps for metric space targets.
\newblock {\em Comm. Anal. Geom.}, 1(3-4):561--659, 1993.

\bibitem{Lem78}
Luc Lemaire.
\newblock Applications harmoniques de surfaces riemanniennes.
\newblock {\em J. Differential Geom.}, 13(1):51--78, 1978.

\bibitem{Lem82}
Luc Lemaire.
\newblock Boundary value problems for harmonic and minimal maps of surfaces
  into manifolds.
\newblock {\em Ann. Scuola Norm. Sup. Pisa Cl. Sci. (4)}, 9(1):91--103, 1982.

\bibitem{LW16-harmonic}
Alexander Lytchak and Stefan Wenger.
\newblock Regularity of harmonic discs in spaces with quadratic isoperimetric
  inequality.
\newblock {\em Calc. Var. Partial Differential Equations}, 55(4):55:98, 2016.

\bibitem{LW15-Plateau}
Alexander Lytchak and Stefan Wenger.
\newblock Area minimizing discs in metric spaces.
\newblock {\em Arch. Ration. Mech. Anal.}, 223(3):1123--1182, 2017.

\bibitem{LW-intrinsic}
Alexander Lytchak and Stefan Wenger.
\newblock Intrinsic structure of minimal discs in metric spaces.
\newblock {\em Geom. Topol.}, 22(1):591--644, 2018.

\bibitem{LWY20}
Alexander Lytchak, Stefan Wenger, and Robert Young.
\newblock Dehn functions and {H}\"{o}lder extensions in asymptotic cones.
\newblock {\em J. Reine Angew. Math.}, 763:79--109, 2020.

\bibitem{MZ10}
Chikako Mese and Patrick~R. Zulkowski.
\newblock The {P}lateau problem in {A}lexandrov spaces.
\newblock {\em J. Differential Geom.}, 85(2):315--356, 2010.

\bibitem{Mor48}
Charles~B. Morrey, Jr.
\newblock The problem of {P}lateau on a {R}iemannian manifold.
\newblock {\em Ann. of Math. (2)}, 49:807--851, 1948.

\bibitem{Nik79}
I.~G. Nikolaev.
\newblock Solution of the {P}lateau problem in spaces of curvature at most
  {$K$}.
\newblock {\em Sibirsk. Mat. Zh.}, 20(2):345--353, 459, 1979.

\bibitem{OvdM14}
Patrick Overath and Heiko von~der Mosel.
\newblock Plateau's problem in {F}insler 3-space.
\newblock {\em Manuscripta Math.}, 143(3-4):273--316, 2014.

\bibitem{Rad30}
Tibor Rad{\'o}.
\newblock On {P}lateau's problem.
\newblock {\em Ann. of Math. (2)}, 31(3):457--469, 1930.

\bibitem{Res97}
Yu.~G. Reshetnyak.
\newblock Sobolev classes of functions with values in a metric space.
\newblock {\em Sibirsk. Mat. Zh.}, 38(3):657--675, iii--iv, 1997.

\bibitem{Res06}
Yu.~G. Reshetnyak.
\newblock On the theory of {S}obolev classes of functions with values in a
  metric space.
\newblock {\em Sibirsk. Mat. Zh.}, 47(1):146--168, 2006.

\bibitem{SU82}
J.~Sacks and K.~Uhlenbeck.
\newblock Minimal immersions of closed {R}iemann surfaces.
\newblock {\em Trans. Amer. Math. Soc.}, 271(2):639--652, 1982.

\bibitem{SY79}
R.~Schoen and Shing~Tung Yau.
\newblock Existence of incompressible minimal surfaces and the topology of
  three-dimensional manifolds with nonnegative scalar curvature.
\newblock {\em Ann. of Math. (2)}, 110(1):127--142, 1979.

\bibitem{Shi39}
Max Shiffman.
\newblock The {P}lateau problem for minimal surfaces of arbitrary topological
  structure.
\newblock {\em Amer. J. Math.}, 61:853--882, 1939.

\bibitem{TT88}
Friedrich Tomi and Anthony~J. Tromba.
\newblock Existence theorems for minimal surfaces of nonzero genus spanning a
  contour.
\newblock {\em Mem. Amer. Math. Soc.}, 71(382):iv+83, 1988.

\bibitem{White86}
Brian White.
\newblock Infima of energy functionals in homotopy classes of mappings.
\newblock {\em J. Differential Geom.}, 23(2):127--142, 1986.

\bibitem{White88}
Brian White.
\newblock Homotopy classes in {S}obolev spaces and the existence of energy
  minimizing maps.
\newblock {\em Acta Math.}, 160(1-2):1--17, 1988.

\end{thebibliography}
\end{document}